\documentclass[arial,10pt]{article}

\usepackage{newalg}
\usepackage{bm}
\usepackage{graphicx}
\usepackage{verbatim} 
 
\input epsf.sty
\input psfig.sty

\def\comment#1{}

\date{}
\begin{document}

\title{Evolving Graph Representation and Visualization}

\author{Anurat Chapanond, Mukkai S. Krishnamoorthy      \\
Department of Computer Science  \\
Rensselaer Polytechnic Institute,       \\
Troy, NY 12180, \\
chapaa@cs.rpi.edu; moorthy@cs.rpi.edu
\and G. M. Prabhu               \\
Department of Computer Science     \\
Iowa State University              \\
Ames, IA 50011                     \\
prabhu@cs.iastate.edu
\and J. Punin                   \\
Software Engineer                  \\
Oracle                             \\
Redwood Shores, CA 94065           \\
puninj@gmail.com}

\maketitle

\section*{Abstract}

The study of evolution of networks has received increased interest with the recent discovery that many real-world networks possess many things in common, in particular the manner of evolution of such networks. By adding a dimension of time to graph analysis, evolving graphs present opportunities and challenges to extract valuable information. This paper introduces the Evolving Graph Markup Language (EGML), an XML application for representing evolving graphs and related results. Along with EGML, a software tool is provided for the study of evolving graphs. New evolving graph drawing techniques based on the force-directed graph layout algorithm are also explored. Our evolving graph techniques reduce vertex movements between graph instances, so that an evolving graph can be viewed with smooth transitions.

\section{Introduction}

Building visualization tools for networks is an active research topic. An effective visualization tool is useful for understanding the mutations of an evolving graph. Picturing data as images (or movies) is more helpful than summary statistics. 

Visualization tools for evolving graphs can be classified into two categories \cite{Moody}. The first category views evolving graphs as several instances of static graphs called `snapshots' \cite{Bearman}. Visualization systems have been built wherein evolving graphs are viewed as multiple snapshots of static graphs. This is done by placing the graphs next to each other so that a user can look at the graphs side by side and see the transitions from one graph to the next. If an evolving graph has a small number of vertices, this approach works well. However, the larger the size of the graph, the more difficult it is for a user to distinguish changes between two graphs.

The second category views an evolving graph as a `movie' \cite{Bender-deMoll},
where graph instances are shown one at a time and the next graph instance is shown by gradually modifying the current instance. 

Ideally, vertices which do not undergo change remain {\it insitu}, while vertices that are affected by changes move to their new positions. The research on evolutionary graphs is still evolving and the focus of this paper is on scalability in the implementation.

This paper is organized as follows. Section 2 consists of a brief literature review of some pertinent papers. Section 3 introduces Evolving Graph Markup Language (EGML) and a software tool for visualization of evolving graphs. The datasets used, namely, Eurovision and US House of Representatives, are also explained. Section 4 describes the design of two new evolving graph layout algorithms – vertex optimization and vector optimization. Section 5 contains experimental results on the two new evolving graph layout algorithms on different types of evolving graphs, and Section 6 consists of the conclusions.

\section{Literature Review}

Software for visualization of large networks has been implemented by Batagelj and Mrvar \cite{Batagelj}. The main contribution of this work is scalability because it is very difficult to implement visualization for graphs with more than 1,000 vertices. It is also important to visualize changes in evolving graphs. Preston and Krishnamoorthy initiated methods for graphical visualization of evolving graphs and provided many visualization methods \cite{Preston}. Some visualization work has been done by Punin et al. who provided a graph representation in XML format \cite{Punin}, and by Fruchterman and Reingold \cite{Fruchterman} who proposed the elegant force-directed placement algorithm. Other researchers have contributed to visualization analysis using graph representation, e.g., Cole et al. \cite{Cole} and Wegman and Marchette \cite{Wegman}. A comprehensive survey of graph drawing and algorithms for visualization of graphs is given in the book by Tollis et al. \cite{Tollis}. Recently Moody et al. focused on visualizing dynamic networks and give an overview of various applications in longitudinal social networks \cite{Moody}.

\section{Markup language and software tool for evolving graphs}

\subsection{Evolving Graph Markup Language}

The {\bf Evolving Graph Markup Language (EGML)} is an extended application of 
the eXtensible Graph Markup and Modeling Language (XGMML) \cite{Punin3}, 
which is an XML application based on the Graph Markup Language (GML) 
used for describing graphs. The purpose of EGML is to enable the exchange of evolving graphs between different authoring and browsing tools. There are also other graph markup languages (NaGML for instance) \cite{Bradley},\cite{XGMML}.

EGML DTD (Document Type Definition) specifications and EGML XSD (XML Schema Definition) specifications are provided in the appendix. EGML DTD was developed in the doctoral dissertation of the first author under the supervision of the second author \cite{PuninPhD}.

\subsection{Software tool for evolving graphs}

A special software tool that reads and writes EGML-formatted files has been built to help study evolving graphs. This tool allows the user to view and manipulate evolving graphs. It also allows the user to view evolving graph transitions where each graph instance is shown with vertex and edge transitions
to the next graph instance.

\begin{figure}[hbtp]
\centering 
\includegraphics[width=4.5in]{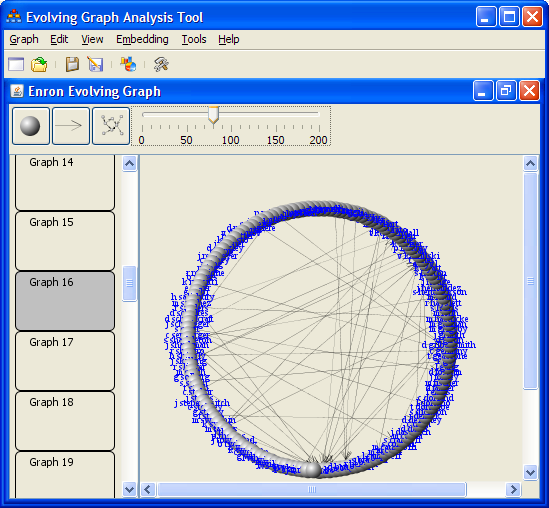}
\caption{The interface for manipulating evolving graphs}
\label{fig:soft}
\end{figure}

Figure \ref{fig:soft} shows the graphical user interface for viewing and manipulating evolving graphs. The evolving graph window, which shows an evolving graph, contains a tool bar on the top, the navigating bar on the left, and the 
evolving graph view on the right. The tool bar contains three buttons and a slider. Each button is used for adding new vertices, edges, or graphs respectively. The slider in the tool bar changes the scale of the evolving graph view. The navigating bar is used for selecting a graph instance in the evolving graph view. The evolving graph view also allows the user to change a vertex's position and vertex and edge properties.

\begin{figure}[hbtp]
\centering 
\includegraphics[width=4.5in]{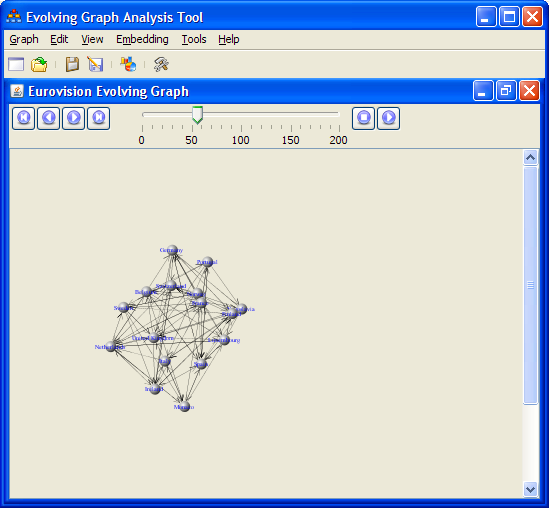}
\caption{The interface for viewing evolving graph transitions}
\label{fig:soft2}
\end{figure}

Figure \ref{fig:soft2} shows the graphical user interface for viewing evolving graph transitions. The transition window contains a tool bar on the top and the transition view. The tool bar contains a group of buttons used to navigate graph instances, a slider for scaling transition view, and two buttons for starting and stopping the transition.

\subsection{Examples of evolving graphs}

\subsubsection{Eurovision dataset}

The Eurovision Song Contest is an annual competition held among European countries. Each country submits a song and casts votes for the other countries' songs, rating which ones are the best. The country with the highest score wins the competition. Broadcast since 1956, this competition is one of the most-watched non-sporting events in the world \cite{Eurovision}.

The rules of the contest have changed over time, though the current scoring
system has been in use since 1975. Each participating country votes
for the ten other countries by assigning the following set of points: 
$ \{ 1, 2, 3, 4, 5, 6, 7, 8, 10, 12 \} $ \cite{Fenn}. 
Watts \cite{Watts} has an interesting view 
on how politics affects vote results.

\begin{figure}[hbtp]
\centering 
\includegraphics[scale=0.25]{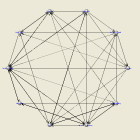}
\includegraphics[scale=0.25]{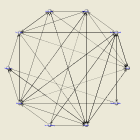}
\includegraphics[scale=0.25]{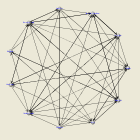}
\includegraphics[scale=0.25]{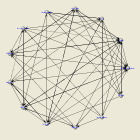}
\includegraphics[scale=0.25]{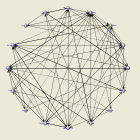}
\includegraphics[scale=0.25]{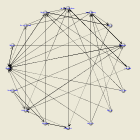}
\includegraphics[scale=0.25]{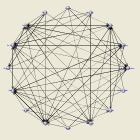}
\includegraphics[scale=0.25]{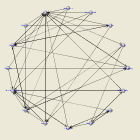}
\includegraphics[scale=0.25]{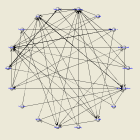}
\includegraphics[scale=0.25]{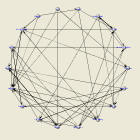}
\\
\includegraphics[scale=0.25]{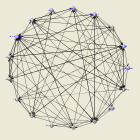}
\includegraphics[scale=0.25]{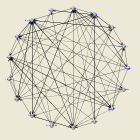}
\includegraphics[scale=0.25]{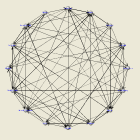}
\includegraphics[scale=0.25]{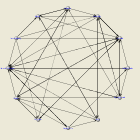}
\includegraphics[scale=0.25]{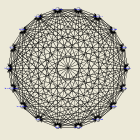}
\includegraphics[scale=0.25]{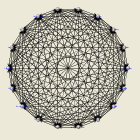}
\includegraphics[scale=0.25]{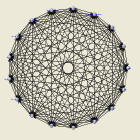}
\includegraphics[scale=0.25]{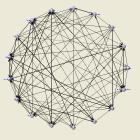}
\includegraphics[scale=0.25]{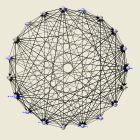}
\includegraphics[scale=0.25]{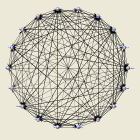}
\\
\includegraphics[scale=0.25]{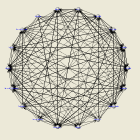}
\includegraphics[scale=0.25]{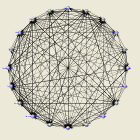}
\includegraphics[scale=0.25]{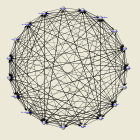}
\includegraphics[scale=0.25]{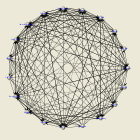}
\includegraphics[scale=0.25]{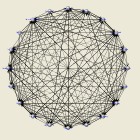}
\includegraphics[scale=0.25]{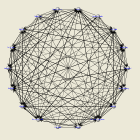}
\includegraphics[scale=0.25]{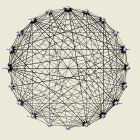}
\includegraphics[scale=0.25]{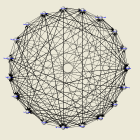}
\includegraphics[scale=0.25]{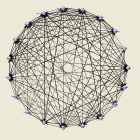}
\includegraphics[scale=0.25]{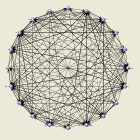}
\\
\includegraphics[scale=0.25]{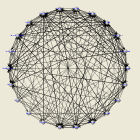}
\includegraphics[scale=0.25]{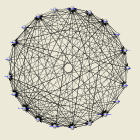}
\includegraphics[scale=0.25]{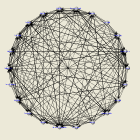}
\includegraphics[scale=0.25]{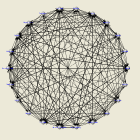}
\includegraphics[scale=0.25]{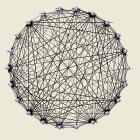}
\includegraphics[scale=0.25]{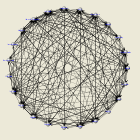}
\includegraphics[scale=0.25]{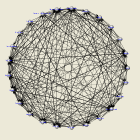}
\includegraphics[scale=0.25]{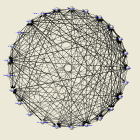}
\includegraphics[scale=0.25]{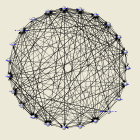}
\includegraphics[scale=0.25]{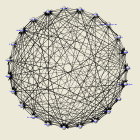}
\\
\includegraphics[scale=0.25]{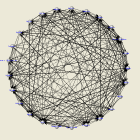}
\includegraphics[scale=0.25]{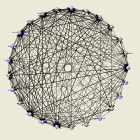}
\includegraphics[scale=0.25]{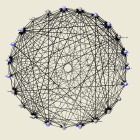}
\includegraphics[scale=0.25]{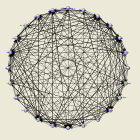}
\includegraphics[scale=0.25]{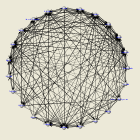}
\includegraphics[scale=0.25]{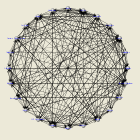}
\includegraphics[scale=0.25]{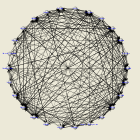}
\\
\caption{Eurovision evolving graph with circular layout}
\label{fig:app_eg_euro}
\end{figure}

The Eurovision evolving graph is constructed from competition data
ranging from 1956 to 2002. Each vertex represents a participating country 
and a directed edge represents a vote a country cast for another country.
There are $ 47 $ graph instances of undirected weighted graphs 
with $ 41 $ distinct countries.
The circular graph layout of Eurovision evolving graph is depicted in 
Figure \ref{fig:app_eg_euro}.

\subsubsection{US House of Representatives dataset}

The US House of Representatives dataset \cite{Rollcall}
contains roll call votes of the representatives from 1989-2004.
The data are stored in text file format and each file contains roll call votes
during a two-year period.
The following shows sample roll call votes data:

\begin{verbatim}
1019990899 0USA     20000BUSH       9999996999169699999996
1011509041 1ALABAMA 20001CALLAHAN  H6616166161166616111166
1011071741 2ALABAMA 20001DICKINSON  6611166166166616666116
1011563241 3ALABAMA 10002BROWDER  GL0000000000000000000000
1011103741 3ALABAMA 10051NICHOLS  BI0000000000000000000000
1011100041 4ALABAMA 10001BEVILL  TOM1166199111616116111166
1011441941 5ALABAMA 10001FLIPPO  RON1166111111116111111116
1011502241 6ALABAMA 10001ERDREICH  B1166111111616116111166
1011541641 7ALABAMA 10001HARRIS  CLA1166111111116116111116
1011406681 1ALASKA  20001YOUNG  DONA6611611169116111966116
1011544061 1ARIZONA 20001RHODES  JOH6616166166166611661116
1011056661 2ARIZONA 10001UDALL  MORR1161111119611111116116
1011445461 3ARIZONA 20001STUMP  BOB 6616166166166616661161
1011542961 4ARIZONA 20001KYL  JON   6616166166166616666161
\end{verbatim}

Each line of text contains vote results for a member of the House.
The vote result, starting at column $ 37 $ (not counting the blank spaces), can be one of the following four categories:

\begin{itemize}
\item $ 0 $: Not present
\item $ 1 $: Agree
\item $ 6 $: Reject
\item $ 9 $: Present, not voting
\end{itemize}

The US house of Representatives evolving graph is constructed with
each vertex representing a member of the House, each edge
representing similar roll call votes between two representatives,
and an edge weight representing the number of similar roll call votes.
Only the vote results $ 1 $ and $ 6 $ are used to construct an edge.

\begin{figure}[hbtp]
\centering 
\includegraphics[scale=0.2]{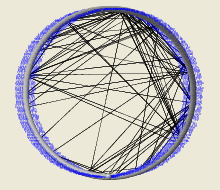}
\includegraphics[scale=0.2]{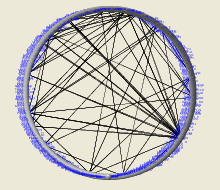}
\includegraphics[scale=0.2]{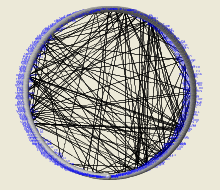}
\includegraphics[scale=0.2]{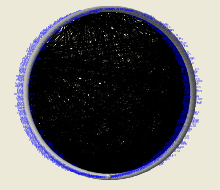}
\includegraphics[scale=0.2]{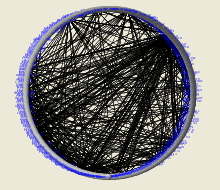}
\includegraphics[scale=0.2]{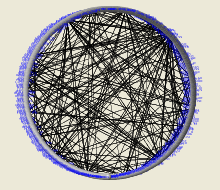}
\includegraphics[scale=0.2]{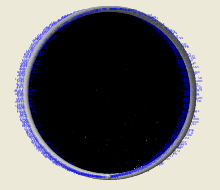}
\includegraphics[scale=0.2]{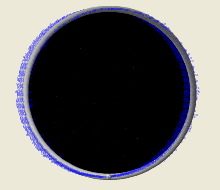}
\caption{US House of Representatives evolving graph with circular layout}
\label{fig:app_eg_us}
\end{figure}

There are $ 8 $ graph instances of undirected weighted graphs with $ 443 $
vertices in each graph instance.
The circular graph layout of the US house of Representatives evolving 
graph is in Figure \ref{fig:app_eg_us}.

\section{Graph Drawing Algorithms}

\subsection{Force directed algorithm }

The force directed algorithm \cite{Fruchterman} is among the most common graph 
layout algorithms. It was designed for undirected graphs and intended
to achieve the even distribution of vertices over the frame, minimizing
edge crossing, making edge length uniform, and reflecting inherent symmetry.

There are two principles for graph drawing in the force directed algorithm.
The first principle is to draw vertices that are connected by edges 
close to each other, and the second principle is that vertices should not be drawn too close to each other.

The force directed algorithm considers vertices behaving as atomic particles, exerting attractive forces and repulsive forces on one another. These attractive and repulsive forces induce movements of the vertices. The attractive force is generated between neighboring vertices and the repulsive force is generated by all vertices. These forces depend on the placement of vertices on a frame. 

The total force obtained by the total sum of vectors generated by repulsive
forces and attractive forces for each vertex is used to determine the new
position in \textsc{Displace-vertex}, together with the temperature value which limits the maximum distance each vertex can move from its current position.

\subsection{Vertex optimization algorithm}

The Vertex optimization algorithm developed by the authors is an algorithm for drawing an evolving graph based on the force directed algorithm. It tries to optimize the vertex position to reduce the amount of movements of each vertex between the transitions. Each instance of evolving graph will be subjected to the force directed algorithm using forces generated by the graph itself and its neighboring graph instances. The vertices of the neighboring graph instances are used to generate forces that draw the same vertices on the current graph instance close to them. The outline of vertex optimization algorithm is the following.

\begin{algorithm}{ Vertex-optimization }{ EG }
\begin{FOR}{ G \in EG }	\\
	\CALL{Vertex-optimized-force-directed}( G )
\end{FOR}
\end{algorithm}

\begin{algorithm}{ Vertex-optimized-force-directed }{ G }
\begin{FOR}{ i \= 1 \TO iterations }	\\
	\begin{FOR}{ u \in V }
		\begin{FOR}{ v \in V }
			\begin{IF}{ u \neq v }
				\CALL{Calculate-repulsive-force}( u, v )
			\end{IF}
		\end{FOR}
	\end{FOR}	\\
	\begin{FOR}{ \mbox{ each edge } ( u, v ) \in E }
		\CALL{Calculate-attractive-force}( u, v )
	\end{FOR}	\\
	\begin{FOR}{ u \in V }
		\begin{FOR}{ v \in \mbox{ neighbor of } G | u = v }
		\CALL{Calculate-attractive-force}( u, v )
		\end{FOR}
	\end{FOR}	\\
	\begin{FOR}{ u \in V }
		\CALL{Displace-vertex}( u, temp )
	\end{FOR}
\end{FOR}
\end{algorithm}

Vertex optimization algorithm runs \textsc{Vertex-optimized-force-directed} for each instance of an evolving graph. It is the same as the force directed algorithm with addition of lines 9 to 11 where it uses vertices from neighboring graph instance to calculate attractive forces between them.

Vertex $ u $ is from graph instance $ G $ and vertex $ v $ is from neighboring graph instance of $ G $. The use of attractive force to pull the same vertex on different graph instances together can be thought of as adding a new edge between them although the function for calculating attractive force may not have to be the same.

{\bf Window Size}

The vertex optimization algorithm uses vertices from neighboring graph instances to determine placement of vertices. As the goal of the algorithm is to find the optimal position of vertices for each graph instance, it is not necessary to use only vertices from neighboring graph instances, which are the graph instances that are next to the considered graph instance. The more number of graph instances used, the better the vertex positions determined. 

Window size is the number of graph instances that the algorithm uses to generate the attractive forces.

\subsection{Vector optimization algorithm}

The vector optimization algorithm, also developed by the authors, addresses the problem of making the transition between graph instances smoother by considering the direction that vertices move from and to. In vertex optimization algorithm, the goal is to minimize the distance between the same vertex on adjacent graph instances. However, in vector optimization algorithm, not only the distance is considered, but also the direction to which the vertex moves. The movement is considered smooth if the direction that the vertex moves to is somewhat the same direction. For example if a vertex moves from left to right but then moves to left, this makes the movement not smooth, since the vertex moves in the opposite direction. 

The idea of vector optimization is to penalize more a vertex which has such movement by adding more force in the opposite direction.

The vector optimization algorithm constructs vectors that represent movement of vertex from vertex position of considered graph instance and its neighboring graph instances. These vectors are used to determine whether a vertex has an unsmooth movement and penalty can be assigned to the vertex.
The outline for the vector optimization algorithm is as follows.

\begin{algorithm}{ Vector-optimization }{ EG }
\begin{FOR}{ G \in EG }	\\
	\CALL{Vector-optimized-force-directed}( G_p, G, G_n )
\end{FOR}
\end{algorithm}

\begin{algorithm}{ Vector-optimized-force-directed }{ G_p, G, G_n }
\begin{FOR}{ i \= 1 \TO iterations }	\\
	\begin{FOR}{ u \in V }
		\begin{FOR}{ v \in V }
			\begin{IF}{ u \neq v }
				\CALL{Calculate-repulsive-force}( u, v )
			\end{IF}
		\end{FOR}
	\end{FOR}	\\
	\begin{FOR}{ \mbox{ each edge } ( u, v ) \in E }
		\CALL{Calculate-attractive-force}( u, v )
	\end{FOR}	\\
	\begin{FOR}{ u \in V }
		u_p \in V_p | u = u_p	\\
		u_n \in V_n | u = u_n	\\
		\nu_1 = \CALL{construct-vector} ( u_p, u )	\\
		\nu_2 = \CALL{construct-vector} ( u, u_n )	\\
		\nu_a = \nu_1 +\nu_2	\\
		\begin{IF}{ |\nu_a| < \mbox{max}( \nu_1, \nu_2 ) }
			\CALL{Calculate-attractive-force}( u, v, penalty )
		\ELSE
			\CALL{Calculate-attractive-force}( u, v )
		\end{IF}
	\end{FOR}	\\
	\begin{FOR}{ u \in V }
		\CALL{Displace-vertex}( u, temp )
	\end{FOR}
\end{FOR}
\end{algorithm}

The Vector optimization algorithm runs 
\textsc{Vector-optimized-force-directed}\\ for each instance of an evolving graph. It is the same as the force directed
algorithm with addition of lines 9 to 17. It has parameters $ G_p $, $ G $,
and $ G_n $ where $ G_p $ denotes a graph $ G_p = ( V_p, E_p ) $ as
the previous graph instance, $ G $ denotes a graph $ G = ( V, E ) $ as
the current graph instance, and $ G_n $ denotes a graph $ G_n = ( V_n, E_n ) $
as the next graph instance.

Vectors $ \nu_1 $, $ \nu_2 $ and $ \nu_a $ are created from the position
of vertex $ u $ on different graph instances. Vector $ \nu_1 $ is created
from the position of vertex $ u_p $, the vertex $ u $ on the previous graph 
instance and vertex $ u $. 
Vector $ \nu_2 $ is created from the position of vertex $ u $ and 
vertex $ u_n $, the vertex $ u $ on the next graph instance.
Vector $ \nu_a $ is the vector addition of vector $ \nu_1 $ and $ \nu_2 $.

The condition for penalizing a vertex with more force to the opposite 
direction is $ \nu_a < \mbox{max}( \nu_1, \nu_2 ) $, the
addition of two vectors has a size smaller than either of $ \nu_1 $ or
$ \nu_2 $. This is because the size of the vector addition will increase
if two vectors are in about the same direction and the size will
decrease if two vectors are not in the same direction.
The penalized force is computed when the condition holds true in
\textsc{Calculate-attractive-force}. By adding a $ penalty $ value,
the force in opposite direction is computed.

\subsection{Measurement for optimization}

The goal for optimization is to place vertices in a way that minimizes their movements. An obvious parameter for optimization is the total distances these vertices move. The total distances can be measured in different way depending on whether the interest is in the overall performance, the specific vertex, or specific graph.

{\bf Total distance measurement ( $ td_{EG} $ ) }

Total distance measurement adds the movements of all vertices during the transition for all graph instances. The total distance measurement of an evolving graph can be computed as follows:

\begin{algorithm}{ Total-distance }{ EG }
td_{EG} \= 0	\\
\begin{FOR}{ i \= 1 \TO p-1 }	\\
	G_i = ( V_i, E_i )	\in EG	\\
	G_{i+1} = ( V_{i+1}, E_{i+1} ) \in EG	\\
	\begin{FOR}{ u_{i,j} \in V_i }
		\begin{FOR}{ u_{i+1,j} \in V_{i+1} }
			td_{EG} \= td_{EG} +d( u_{i,j}, u_{i+1,j} )
		\end{FOR}
	\end{FOR}	\\
\end{FOR}
\RETURN td_{EG}
\end{algorithm}

{\bf Total distance per graph measurement ( $ td_{G} $ ) }

Total distance per graph measurement sums all the movements of vertices during the transition from graph instance $ G $ to the next graph instance. It can be computed as follows:

\begin{algorithm}{ Total-distance-per-graph }{ G, EG }
td_{G} \= 0	\\
G = ( V, E )	\in EG	\\
G_{n} = ( V_{n}, E_{n} ) \in EG | G_n \mbox{ is the next graph instance from } G	\\
\begin{FOR}{ u_{j} \in V }
	\begin{FOR}{ u_{n,j} \in V_{n} }
		td_{G} \= td_{G} +d( u_{j}, u_{n,j} )
	\end{FOR}
\end{FOR}	\\
\RETURN td_{G}
\end{algorithm}

{\bf Total distance per vertex measurement ( $ td_{v} $ ) }

Total distance per vertex measurement sums all the movements of a vertex during the transition for all graph instances. It can be computed as follows:

\begin{algorithm}{ Total-distance-per-vertex }{ v, EG }
td_{v} \= 0	\\
\begin{FOR}{ i \= 1 \TO p-1 }	\\
	G_i = ( V_i, E_i )	\in EG	\\
	G_{i+1} = ( V_{i+1}, E_{i+1} ) \in EG	\\
	\begin{FOR}{ u_{i,j} \in V_i }
		\begin{FOR}{ u_{i+1,j} \in V_{i+1} }
			td_{v} \= td_{v} +d( u_{i,j}, u_{i+1,j} )
		\end{FOR}
	\end{FOR}	\\
\end{FOR}
\RETURN td_{v}
\end{algorithm}

\section{Experimental results}

The Vertex optimization and Vector optimization algorithms were run on different types of evolving graphs. First, vertex optimization was run on generated evolving graphs which are randomly generated evolving graphs with different degree distribution. Second, vertex optimization was run on five evolving graphs. These evolving graphs were created specifically to understand the effect of the algorithm on evolving graphs. Third, the comparison between vertex optimization and vector optimization was done with evolving graphs. Fourth, vertex optimization algorithm was run on real-world evolving graphs such as Eurovision evolving graphs and US House of Representatives evolving graphs.

\subsection{Vertex optimization on generated evolving graphs}

In this experiment, the vertex optimization algorithm is tested on different types of generated evolving graphs. The effect of window size and performance is also tested.

Three types of generated evolving graphs are used, namely, random evolving graphs, exponential evolving graphs and scale-free evolving graphs. These evolving graphs have their graph instances with different degree distribution. Random evolving graphs are constructed such that their graph instances have Poisson degree distribution. Exponential evolving graphs have graph instances with exponential degree distribution. And scale-free evolving graphs have graph instances with power-law degree distribution.

All evolving graphs are constructed with comparable size. There are 20 graph instances in total where the first graph instance for random evolving graphs has 50 vertices and 20 edges and 5 more edges are added to each graph instance to create the next graph instance.

For exponential and scale-free evolving graphs, the first graph instance has 31 vertices and one more vertex is added for the next graph instance. Therefore the last graph instance will also have 50 vertices which is the same as number of vertices for random evolving graphs.

These generated evolving graphs are tested with the vertex optimization algorithm with different window size and the total distance measurement ( $ td_{EG} $ ) is measured.

\begin{figure}[hbtp]
\centering 
\includegraphics[width=4.5in]{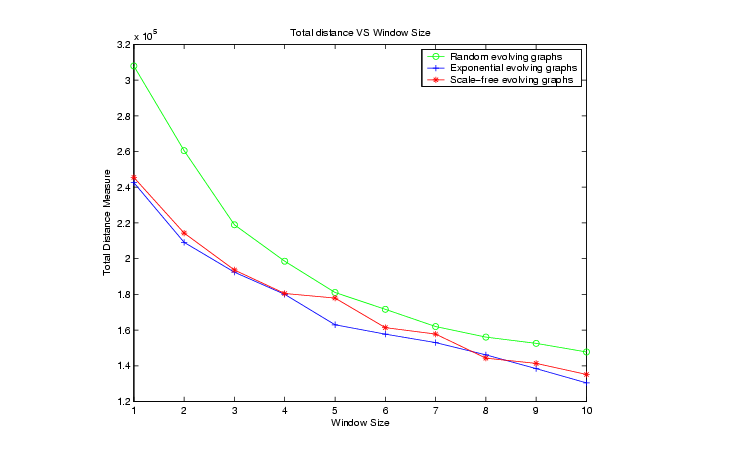}
\caption{The total distance measure with different window size on different types of evolving graphs}
\label{fig:vs_1_distance}
\end{figure}

Figure \ref{fig:vs_1_distance} shows the plot of total distance measure versus varying window size after applying vertex optimization algorithm on the evolving graphs. The different types of evolving graphs are shown as different plots. All types of evolving graphs have the same effect as the total distance decreases when the window size increases. This means with higher window size value the movement distance of vertices become shorter.

Window size indicates the number of graph instances that are used to compute vertex placement and therefore the larger the window size the more the information available for determining the vertex placement.

The random evolving graphs also have higher total distance value than the other two evolving graphs. This is predictable considering the number of vertices the random evolving graphs have. Random evolving graphs have fixed number of vertices for all graph instances and have more vertices than the other two evolving graphs. These additional vertices must contribute to the higher total distance measure.

The exponential evolving graphs and scale-free evolving graphs have the 
same number of vertices and the total distance measure is about the same.

\begin{figure}[hbtp]
\centering 
\includegraphics[width=4.5in]{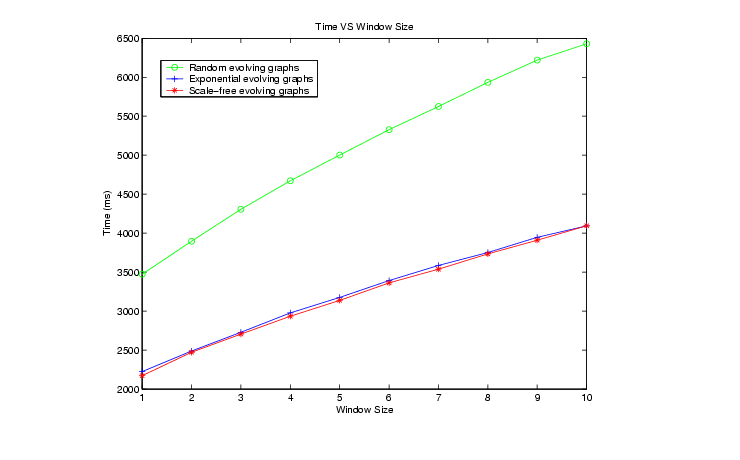}
\caption{The computation time with different window size on different types of evolving graphs}
\label{fig:vs_1_time}
\end{figure}

Figure \ref{fig:vs_1_time} shows the plot of computation time versus varying window size. The plot shows that large window size requires more computation time than small window size. The computation times for random evolving graphs are higher than the other two evolving graphs, also because of the number of vertices as in Figure \ref{fig:vs_1_distance}. These two plots show the trade-off between the optimized vertex position and computation time when changing window size. With larger window size, a better vertex position can be obtained at the expense of more computation time.

\begin{figure}[hbtp]
\centering 
\includegraphics[width=4.5in]{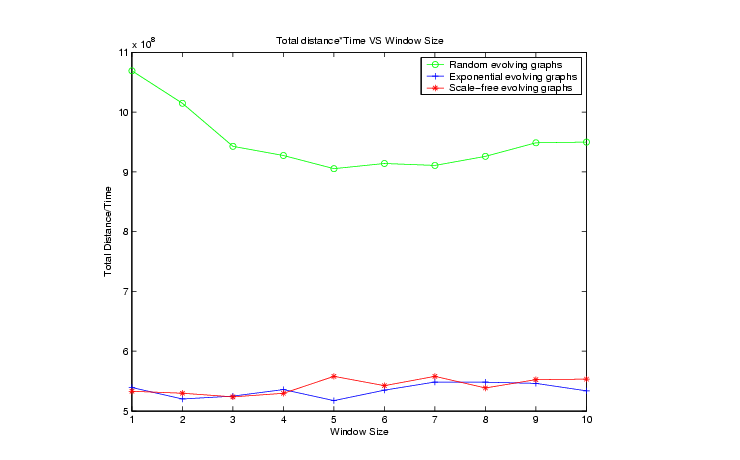}
\caption{The total (distance * computation) time with different window size on different types of evolving graphs}
\label{fig:vs_1_perf}
\end{figure}

Figure \ref{fig:vs_1_perf} shows the plot of total distance measure times the computation time versus varying window size. The value represents the performance trade-off between the optimal position and computation time. 
The lower value for total distance and computation time is preferable; therefore, the lowest point of the plot will represent the optimal window size.
For random evolving graphs, the plot clearly shows the optimal window size at 5. The other plots have no clear optimal window size as the plots are almost straight lines.

This suggests that for random evolving graphs the best window size is around 5 but for real-world evolving graphs which have exponential or power-law distribution, choosing the window size depends more on the trade-off between the cost and computation time i.e., if we are more concerned about the cost, the larger window size is better, and if computation time is more important then lower window size is better.

\subsection{Vertex optimization on evolving graphs}

More examples of evolving graphs are constructed to understand the vertex 
optimization algorithm. All the evolving graphs have twenty graph instances and are described below.

{\bf Example: Evolving graph 1 }

Each instance is a line graph with 53 vertices and edges 
$ ( i ,  i+1  ) | \forall i=1..52 $. 

\begin{figure}[hbtp]
\centering 
\includegraphics[scale=0.5]{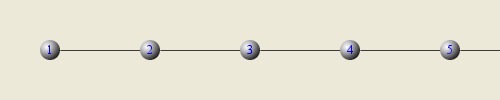}
\caption{Portion of evolving graph 1}
\label{fig:vs_2_ex1}
\end{figure}

Evolving graph 1 is a line graph with no change between graph instances. All the vertices remain unchanged in position in different graph instances. Figure \ref{fig:vs_2_ex1} shows one graph instance of evolving graph 1.Only the first five vertices of the graph instance are shown.

{\bf Example: Evolving graph 2 }

The odd instances are line graphs with 53 vertices and edges 
$ ( i ,  i+1  ) | \forall i=1..52 $. The even instances are graphs with 53 vertices and no edges.

\begin{figure}[hbtp]
\centering 
\includegraphics[scale=0.5]{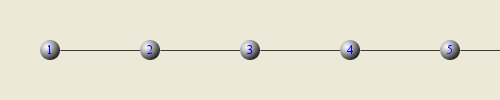}
\\
\includegraphics[scale=0.5]{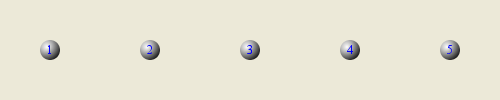}
\caption{Portion of evolving graph 2}
\label{fig:vs_2_ex2}
\end{figure}

Evolving graph 2 has graph instances which switch between a line graph and a graph with no edges. All vertices remain unchanged in position in different graph instances. Figure \ref{fig:vs_2_ex2} shows two graph instances in evolving graph 2. The odd graph instance is on top and the even graph instance is on the bottom in the figure. Only the first five vertices of each graph instance are shown.

{\bf Example: Evolving graph 3 }

The first 15 instances are line graphs with 53 vertices and edges 
$ ( i ,  i+1  ) | \forall i=1..52 $. The last 5 instances are line graphs
with 53 vertices and edges $ ( i, (i+2)\%53 ) | \forall i=1..51,53 $.

\begin{figure}[hbtp]
\centering 
\includegraphics[scale=0.5]{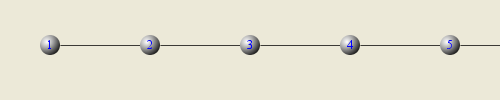}
\\
\includegraphics[scale=0.5]{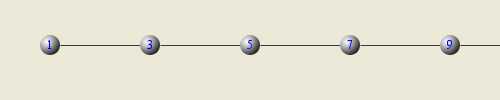}
\caption{Portion of evolving graph 3}
\label{fig:vs_2_ex3}
\end{figure}

Evolving graph 3 has two sets of line graphs. The first set has 15 graph instances and the other set has 5 graph instances. The setup for the second set of graph instances has different vertices in the same position as vertices in the first set of graph instances, e.g., vertex 1 is in the same position for both sets, vertex 2 from the first set is in the same position as vertex 3 from the second set and so on. Figure \ref{fig:vs_2_ex3} shows two graph instances  in evolving graph 3. The first 15 graph instances are the same as the top figure and the last 5 graph instances are the same as the bottom figure. Only the first five vertices of each graph instance are shown.

{\bf Example: Evolving graph 4 }

Each graph instance has 53 vertices. The $ k^{th} $ graph instance contains edges $ ( i, (i+k)\%53 ) | \forall i=1..(53-k) $.

\begin{figure}[hbtp]
\centering 
\includegraphics[scale=0.5]{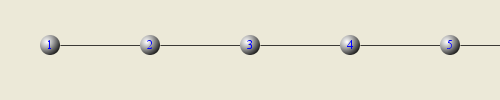}
\\
\includegraphics[scale=0.5]{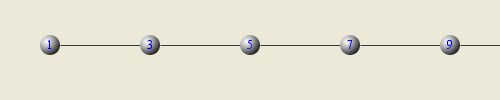}
\\
\includegraphics[scale=0.5]{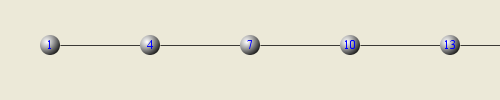}
\\
\includegraphics[scale=0.5]{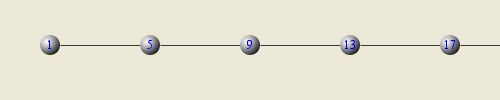}
\\
\includegraphics[scale=0.5]{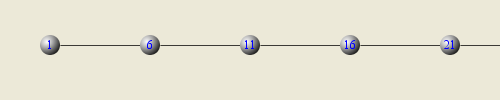}
\\
\caption{Portion of evolving graph 4}
\label{fig:vs_2_ex4}
\end{figure}

Evolving graph 4 has twenty different sets of line graphs. The first line graph has vertices in the order 1, 2, 3, and so on. The second line graph has vertices in the order 1, 3, 5, and so on. Different vertices from different graphs are in the same position, i.e., vertex 1 is in the same position for all graph instances, vertex 2 from the first graph instance, vertex 3 from the second  graph instance, vertex 4 from the third graph instance are in the position. Figure \ref{fig:vs_2_ex4} shows the first five graph instances in evolving graph 4. Only the first five vertices are shown in each graph instance.

{\bf Example: Evolving graph 5 }

Each graph instance is a circular graph with 53 vertices and edges
$ ( i, i+1 ) | \forall i=1..53 $ and extra edges 
$ ( 40, j ) | \forall j=1..20 $ but the edge $ ( 40, k ) $ for the
$ k^{th} $ graph instance is absent.

\begin{figure}[hbtp]
\centering 
\includegraphics[scale=0.5]{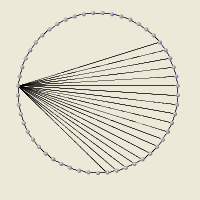}
\includegraphics[scale=0.5]{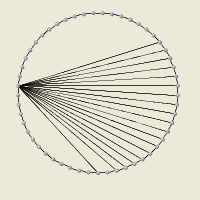}
\includegraphics[scale=0.5]{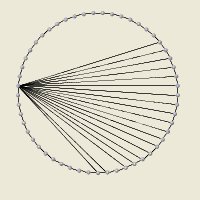}
\includegraphics[scale=0.5]{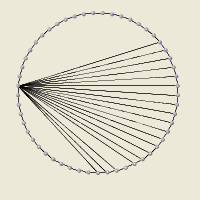}
\caption{Portion of evolving graph 5}
\label{fig:vs_2_ex5}
\end{figure}

Evolving graph 5 has twenty circular graph instances with additional edges from vertex 40 to a set of vertices 1 to 20. One additional edge is missing for each graph instance. All the same vertices in different graph instances are in the same position. Figure \ref{fig:vs_2_ex5} shows the first four graph instances in evolving graph 5.

{\bf Experimental results}

Each evolving graph is run on the vertex optimization algorithm using different window sizes from 1 to 5. The results are quantified by distance measure and plotted in Figure \ref{fig:vs_2_ex1_eg} to \ref{fig:vs_2_ex5_v} in Appendix A.

\subsection{Vertex optimization Vs. Vector optimization algorithm}

The five evolving graphs have been used to compare the two evolving graph layout algorithms.

\begin{figure}[hbtp]
\centering 
\includegraphics[width=4.5in]{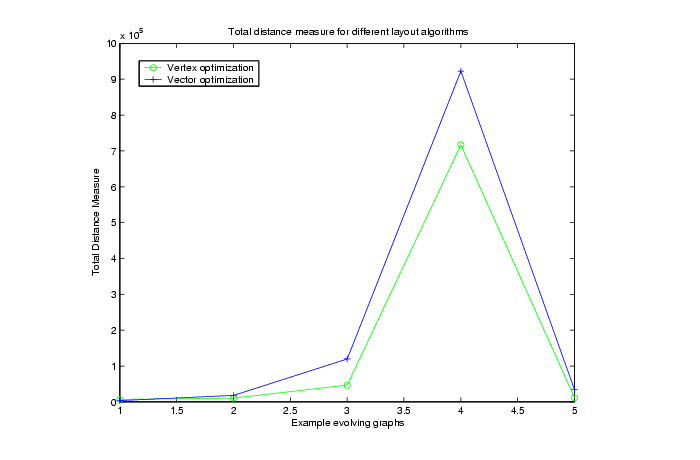}
\caption{The total distance measure with different layout algorithms}
\label{fig:vs_3_n_vs_v_eg}
\end{figure}

Figure \ref{fig:vs_3_n_vs_v_eg} shows the comparison of total distance measure between vertex optimization and vector optimization algorithms. The window size is set to five for both algorithms. The plot shows that vector optimization algorithm does not improve over the vertex optimization algorithm.

\subsection{Vertex optimization on Eurovision and US House of Representatives evolving graphs}

{\bf Eurovision evolving graph}

The circular graph layout of Eurovision evolving graph is in Figure \ref{fig:app_eg_euro}. This layout gives a general view of each graph instance in the Eurovision evolving graph. The first 14 graph instances have smaller number of vertices and edges than the other instances. This is attributed to the change in competition format which results in more consistent number of participating countries and scores each country gives to the others.

\begin{figure}[hbtp]
\centering 
\includegraphics[scale=0.2]{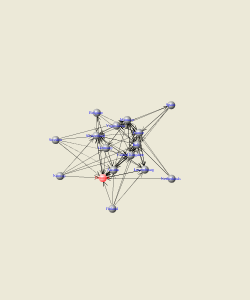}
\includegraphics[scale=0.2]{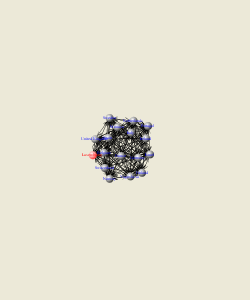}
\includegraphics[scale=0.2]{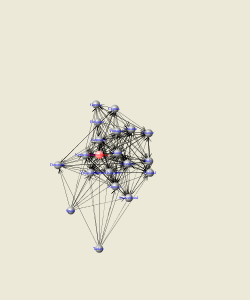}
\includegraphics[scale=0.2]{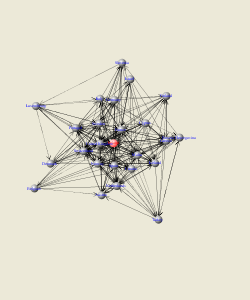}
\includegraphics[scale=0.2]{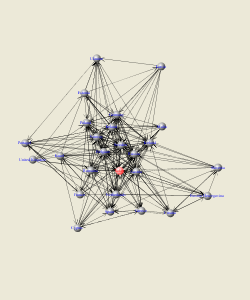}
\\
\caption{Eurovision evolving graph with vertex optimization layout}
\label{fig:app_eg_euro_opt}
\end{figure}

The graph layout of Eurovision evolving graph using vertex optimization algorithm is in Figure \ref{fig:app_eg_euro_opt}. The winning country in each year is represented by the red vertex in the drawing. 
These winning countries, in general, are placed not far from the center of the graph. Although we expect them to be exactly at the center of the graph, that is not the case, which is due to several reasons. First, the competition should be modeled as a weighted directed graph, but vertex optimization does not take edge weight into account. Second, vertex optimization has window size effect which pulls the same vertex in different graph instance together.

The total distance measure for Eurovision evolving graph is as follows:

\begin{figure}[hbtp]
\centering 
\includegraphics[width=4.5in]{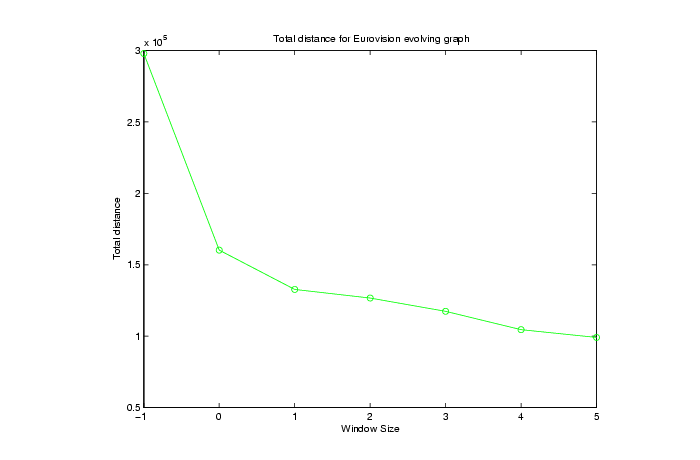}
\caption{The total distance measure with different window sizes on Eurovision evolving graph}
\label{fig:vs_4_euro}
\end{figure}

Figure \ref{fig:vs_4_euro} shows the plot of total distance measure versus varying window sizes on the Eurovision evolving graph. Another two values are added to the plot at window sizes of -1 and 0. The first value describes the total distance when the vertices on Eurovision evolving graph are placed randomly. The second value describes the total distance when the algorithm does not use any information from neighboring graph instances. It is equivalent to running a force directed algorithm on each graph instance separately. Overall, the vertex optimization helps reduce the movement of vertices for Eurovision graphs, the larger the window size, the smaller the total distance measure.

\begin{figure}[hbtp]
\centering 
\includegraphics[width=4.5in]{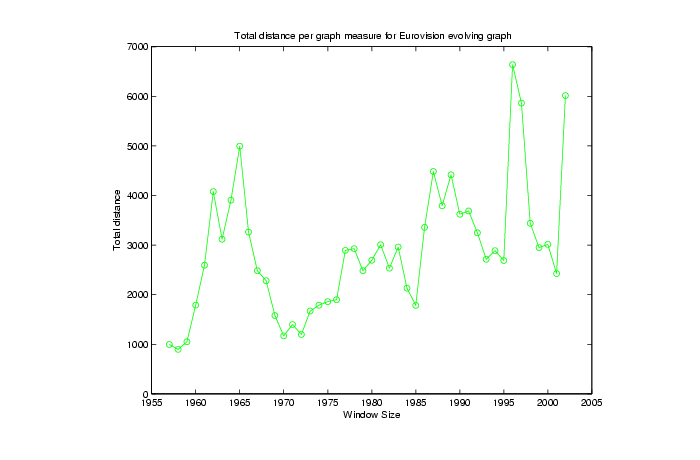}
\caption{The total distance measure per graph on Eurovision evolving graph}
\label{fig:vs_4_euro_g}
\end{figure}

Figure \ref{fig:vs_4_euro_g} shows the plot of total distance measure per graph on Eurovision evolving graph. There are spikes in the plot which is consistent with graph difference.

\begin{table}
\begin{center}
\begin{tabular}{|l|c|c|c|}
\hline
Country 		& Total distance & No.	& Avg total distance \\
\hline
Italy			& 2215		& 31	& 71.45	\\
United Kingdom	& 3390		& 44	& 77.05 \\
Germany			& 3745		& 44	& 85.12 \\
Switzerland		& 3327		& 39	& 85.30 \\
France			& 3877		& 42	& 92.32 \\
Sweden			& 3614		& 39	& 92.67 \\
Romania			& 95.71		&  1	& 95.71 \\
Monaco			& 1929		& 20	& 96.47 \\
Denmark			& 2859		& 27	& 105.89 \\
Netherlands		& 4133		& 38	& 108.77 \\
Luxembourg		& 3781		& 34	& 111.21 \\
Estonia			& 798		&  7	& 114.02 \\
Yugoslavia		& 2813		& 24	& 117.24 \\
Lithuania		& 117		&  1	& 117.99 \\
Ireland			& 4062		& 34	& 119.48 \\
Belgium			& 4838		& 40	& 120.95 \\
Austria			& 4273		& 35	& 122.10 \\
Spain			& 5431		& 42	& 129.31 \\
Norway			& 5338		& 39	& 136.87 \\
Finland			& 4386		& 31	& 141.48 \\
Israel			& 3017		& 21	& 143.70 \\
Malta			& 1942		& 13	& 149.44 \\
Russia			& 657		&  4	& 164.48 \\
Slovenia		& 1065		&  6	& 177.51 \\
Greece			& 3212		& 18	& 178.49 \\
Croatia			& 1930		& 10	& 193.06 \\
Cyprus			& 3610		& 18	& 200.59 \\
Turkey			& 4262		& 21	& 202.96 \\
Iceland			& 2773		& 13	& 213.36 \\
Poland			& 1080		&  5	& 216.02 \\
Portugal		& 7370		& 33	& 223.34 \\
Bosnia Herzegovina	& 1431	&  6	& 238.59 \\
Latvia			& 898		&  3	& 299.65 \\
Hungary			& 811		&  2	& 405.59 \\
FYR Macedonia	& 0			&  0	& 0	\\
Morocco			& 0			&  0	& 0 \\
Slovakia		& 0			&  0	& 0 \\
Ukraine			& 0			&  0	& 0 \\
\hline
\end{tabular}
\caption[Total distance per vertex for each country]
{Total distance per vertex for each country sorted by average total distance}
\label{table:vs_4_euro_v}
\end{center}
\end{table}

Table \ref{table:vs_4_euro_v} shows the total distance per vertex for each country of Eurovision evolving graph. Since some countries do not participate every year, the average total distances per vertex are calculated. The list of countries is sorted by this average. Italy, United Kingdom, Germany, Switzerland, France and Sweden are the six smallest movement countries. Hungary, Latvia, Bosnia Herzegovina, Portugal and Poland are the largest movement countries. FYR Macedonia, Morocco, Slovakia, and Ukraine are countries that never participate in consecutive year so the total distance are not measured.

Countries with small movements are consistent with countries that form a cluster. The total distance per vertex could be another way to determine a cluster in evolving graphs.

{\bf US House of Representatives evolving graph}

The circular graph layout of US House of Representatives evolving graph is in 
Figure \ref{fig:app_eg_us}. This layout shows the density of edges in each graph instance. The $ 4^{th} $, $ 7^{th} $, and $ 8^{th} $ graph instances have high edge density even though edges have been filtered by $ 90\% $ weight threshold.

\begin{figure}[hbtp]
\centering 
\includegraphics[scale=0.25]{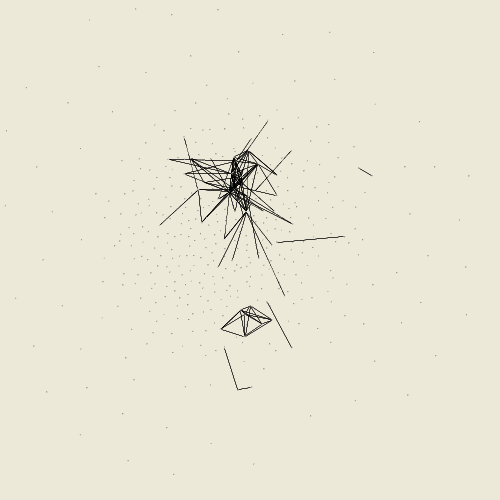}
\includegraphics[scale=0.25]{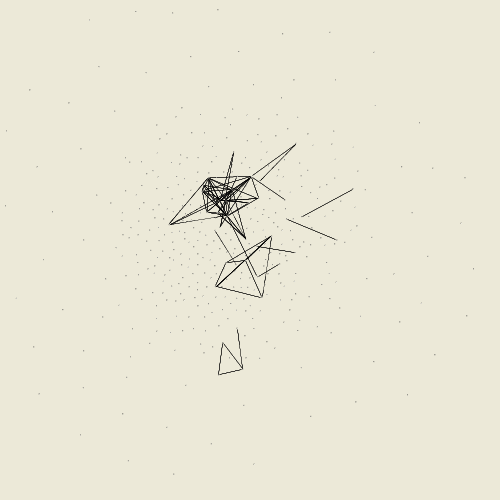}
\\
\includegraphics[scale=0.25]{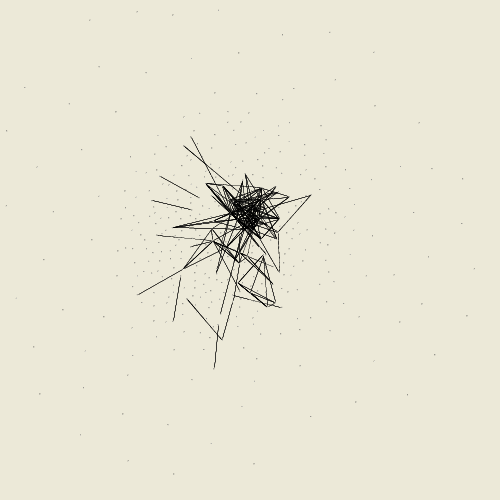}
\includegraphics[scale=0.25]{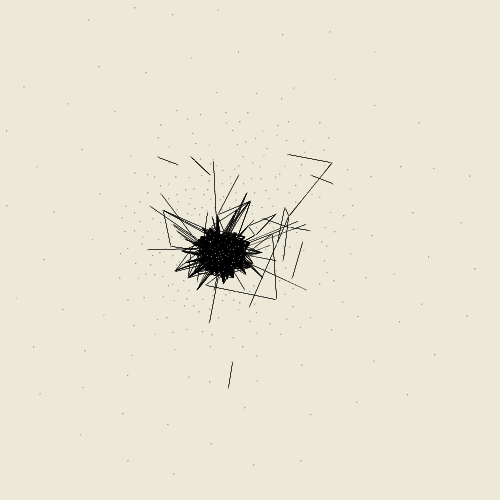}
\\
\includegraphics[scale=0.25]{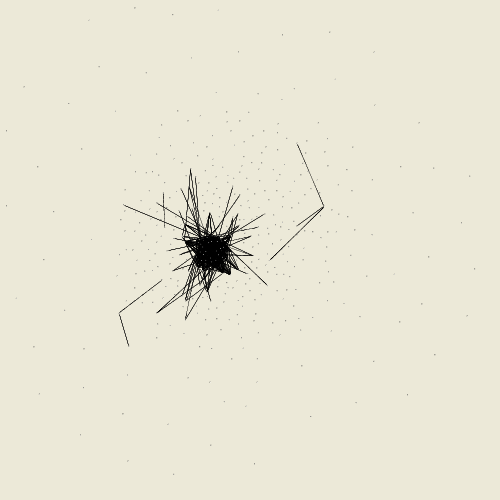}
\includegraphics[scale=0.25]{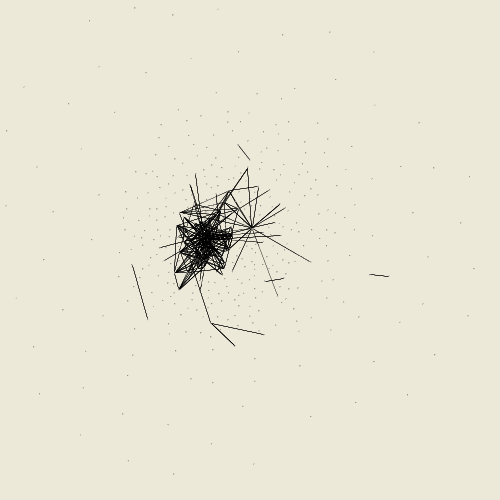}
\\
\includegraphics[scale=0.25]{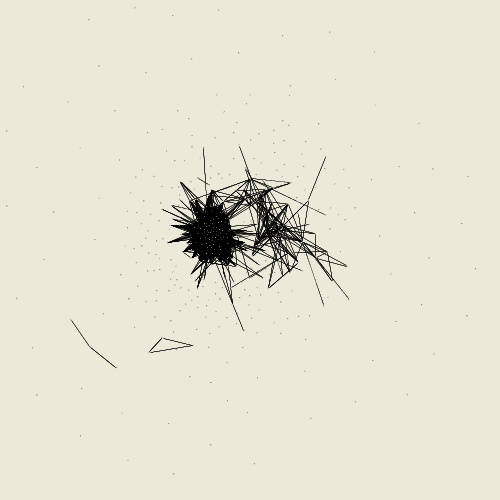}
\includegraphics[scale=0.25]{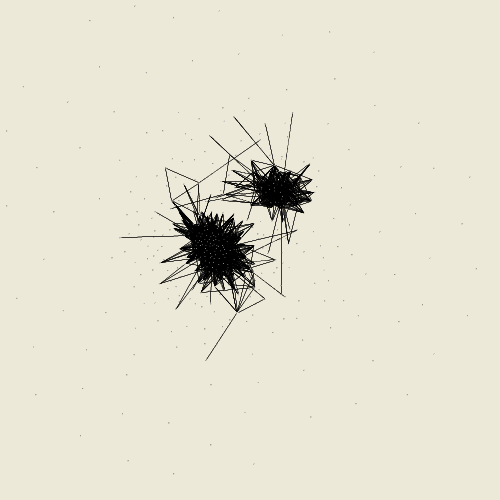}
\\
\caption{US House of Representatives evolving graph with vertex optimization layout}
\label{fig:app_eg_us_opt}
\end{figure}

The graph layout of US House of Representatives evolving graph using vertex optimization algorithm is in Figure \ref{fig:app_eg_us_opt}. Due to very large layout, vertices in the drawing appear as small dots. A giant component can be seen in all graph instances. We expect to see two giant components representing two parties and the layouts show one bigger giant component and a smaller component. The reason the second component cannot be seen clearly is because of the weight threshold. The graph instances in which the second component cannot be seen clearly just do not have enough edges for the second component to form. The graph instances that have high density edges clearly shows two giant components.

The total distance measure for Eurovision evolving graph is as follows:

\begin{figure}[hbtp]
\centering 
\includegraphics[width=4.5in]{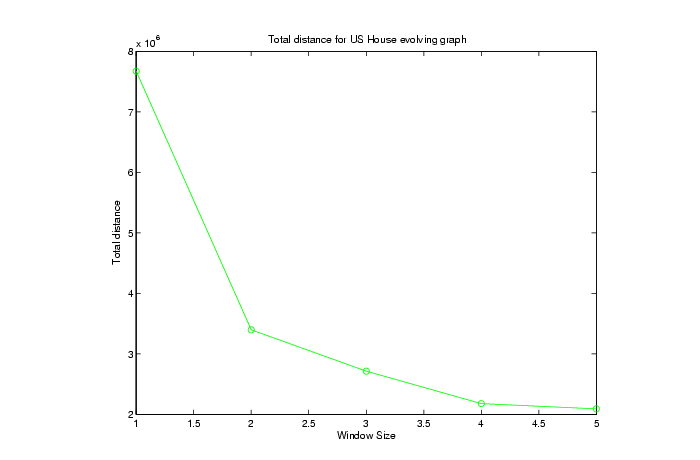}
\caption{The total distance measure on US House of Representatives evolving graph}
\label{fig:vs_4_us_eg}
\end{figure}

Figure \ref{fig:vs_4_us_eg} shows the total distance measure on US House of
Representatives evolving graph with different window sizes. The total distance drops, as the window size increases.

\begin{figure}[hbtp]
\centering 
\includegraphics[width=4.5in]{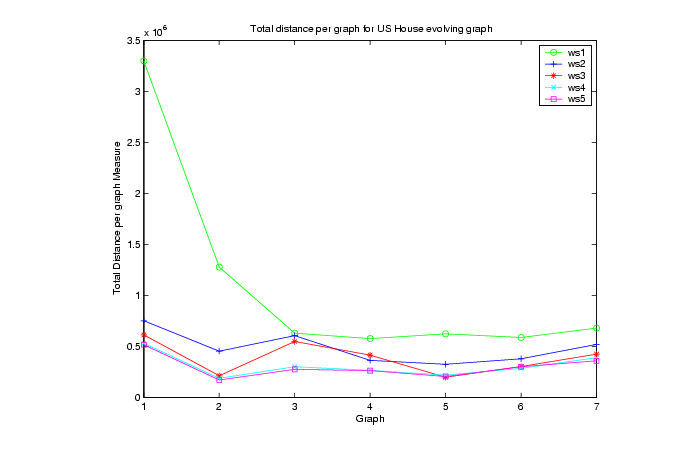}
\caption{The total distance measure per graph on US House of Representatives evolving graph}
\label{fig:vs_4_us_g}
\end{figure}

Figure \ref{fig:vs_4_us_g} shows the total distance per graph measure on US House of Representatives evolving graph. For window sizes higher than one, the total distance measure per graph are quite steady through all graph instances, suggesting that there is no significant difference in graph configuration.

\begin{figure}[hbtp]
\centering 
\includegraphics[width=4.5in]{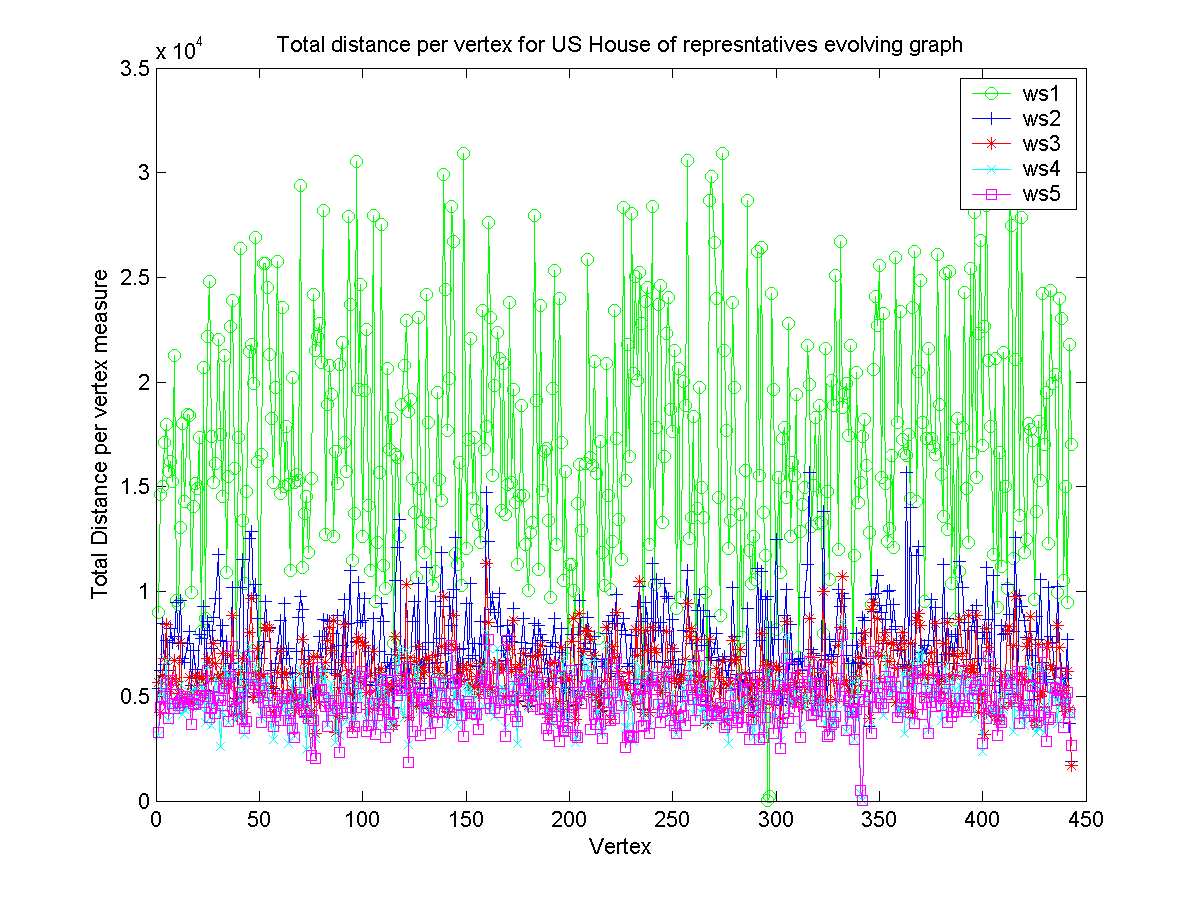}
\caption{The total distance measure per vertex on US House of Representatives evolving graph}
\label{fig:vs_4_us_v}
\end{figure}

Figure \ref{fig:vs_4_us_v} shows the total distance per vertex measure
on US House of Representatives evolving graph. The average total distance per vertex is also calculated. The average value is in the range 265 to 1138. George Bush is in the $ 42^{th} $ place with average value of 470.

\section{Conclusions}

The goal for an evolving graph layout algorithm in the current work is to automatically draw an evolving graph in two dimensions in a manner that augments intuitive understanding of an evolving graph. This involves drawing graph instances to reflect their structure and also providing smooth transitions between graph instances.

Two new evolving graph layout algorithms are introduced and analyzed as part of a research study on visualization tools for evolving graphs. The first algorithm, the vertex optimization algorithm, minimizes the distance between the same vertex in neighboring graph instances by extending the force directed graph layout algorithm. It incorporates attractive forces between the same vertices on neighboring graph instances with the existing forces in force directed algorithm. The second algorithm, the vector optimization algorithm, addresses the smooth movement of vertices between transitions between graph instances. Forces are increased as a penalty for vertices which unsmooth the transitions.

Three measurements are introduced for capturing vertex movements. The total distance explains the overall performance of a layout algorithm. The total distance per graph explains the movement between graph instances as well as the difference between their graph instances' structure. The total distance per vertex explains the impact of different graph instances' structure upon a vertex.

Different types of evolving graphs are tested using vertex optimization algorithm. The experiments show that the difference in degree distribution does not affect the total distance as much as the difference in graph structure. 
Window size, which indicates the amount of information a certain graph instance can gather from its neighboring graph instances, has a significant impact on the total distance measure. However, this comes with a trade-off with computation time. A good window size for small to medium size evolving graphs is five.

Although vector optimization algorithm looks promising to provide a smoother transition between graph instances, the experiments show that there is no improvement over vertex optimization in terms of the total distance measure.

Experiments on real data evolving graphs also show that vertex optimization algorithm can minimize the total movements for vertices in these evolving graphs, thereby giving a smoother transition. The three total distance measurements can be used together with other analysis tools such as clustering, or centrality metrics.

\section*{Appendix A}

{\bf Figures 21 to 36}

\begin{figure}[hbtp]
\centering 
\includegraphics[width=4.5in]{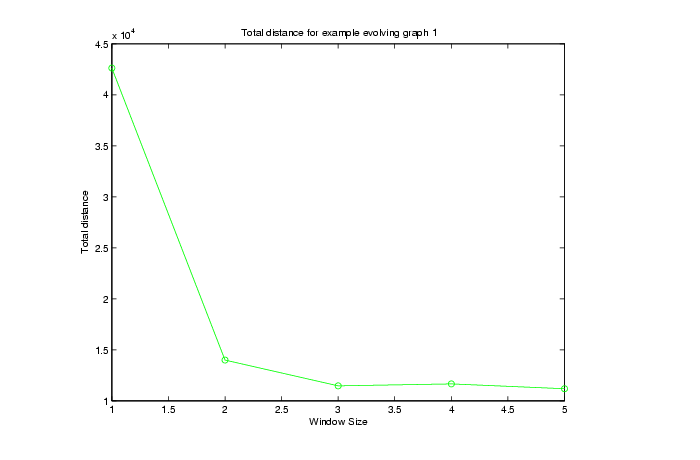}
\caption{The total distance measure with different window sizes on evolving graph 1}
\label{fig:vs_2_ex1_eg}
\end{figure}

\begin{figure}[hbtp]
\centering 
\includegraphics[width=4.5in]{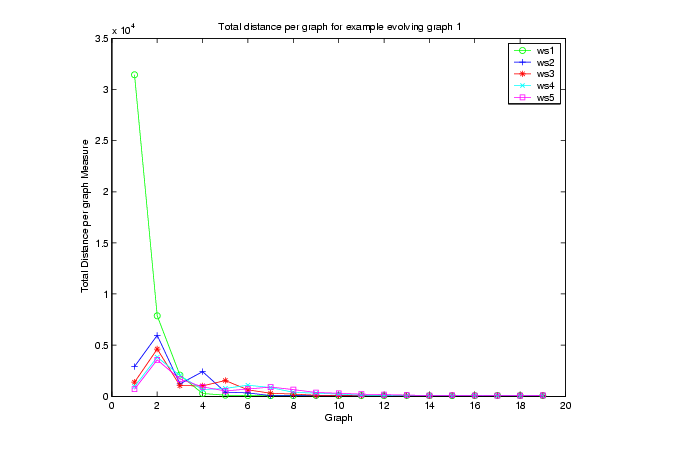}
\caption{The total distance per graph measure with different window sizes on evolving graph 1}
\label{fig:vs_2_ex1_g}
\end{figure}

\begin{figure}[hbtp]
\centering 
\includegraphics[width=4.5in]{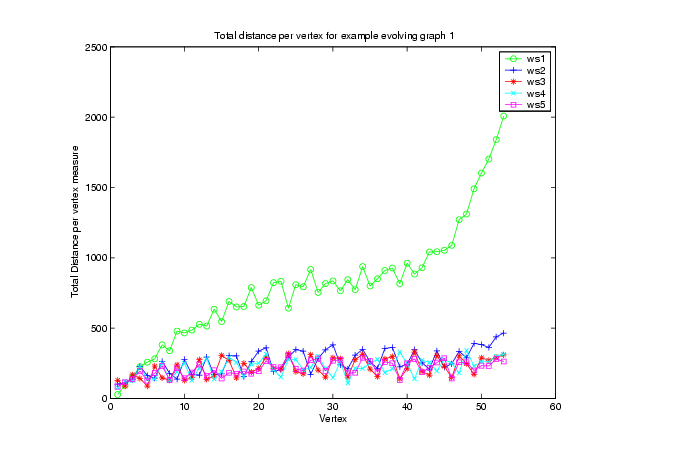}
\caption{The total distance per vertex measure with different window sizes on evolving graph 1}
\label{fig:vs_2_ex1_v}
\end{figure}

\begin{figure}[hbtp]
\centering 
\includegraphics[width=4.5in]{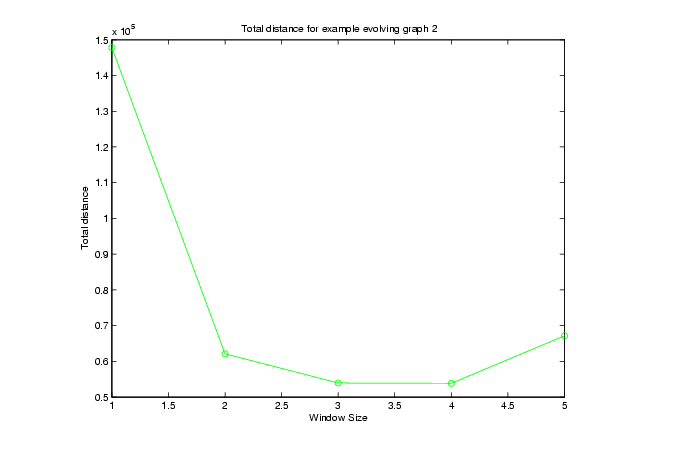}
\caption{The total distance measure with different window sizes on evolving graph 2}
\label{fig:vs_2_ex2_eg}
\end{figure}

\begin{figure}[hbtp]
\centering 
\includegraphics[width=4.5in]{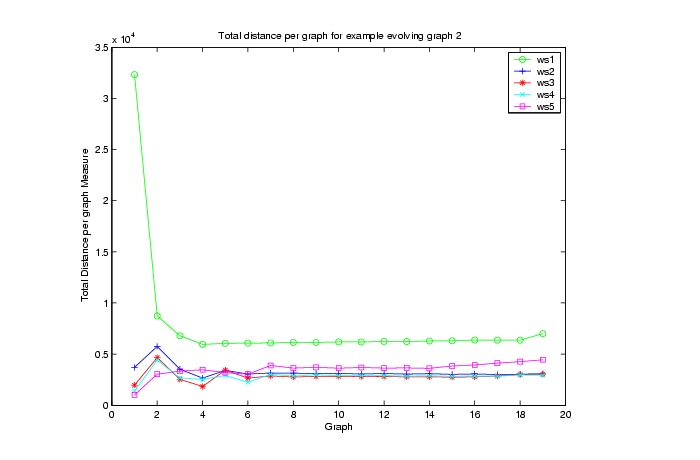}
\caption{The total distance per graph measure with different window sizes on evolving graph 2}
\label{fig:vs_2_ex2_g}
\end{figure}

\begin{figure}[hbtp]
\centering 
\includegraphics[width=4.5in]{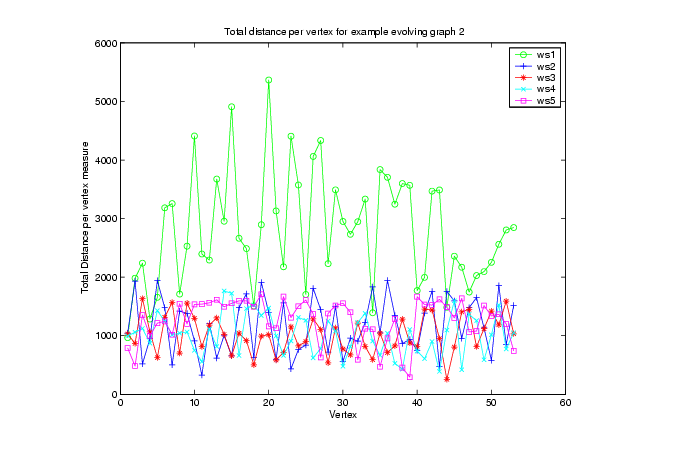}
\caption{The total distance per vertex measure with different window sizes on evolving graph 2}
\label{fig:vs_2_ex2_v}
\end{figure}

\begin{figure}[hbtp]
\centering 
\includegraphics[width=4.5in]{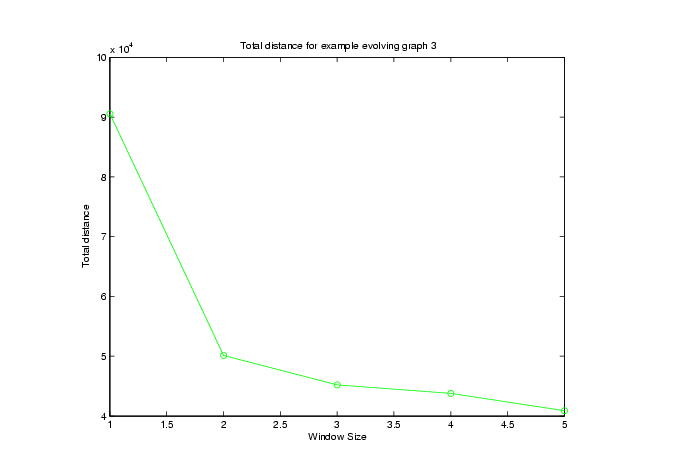}
\caption{The total distance measure with different window sizes on evolving graph 3}
\label{fig:vs_2_ex3_eg}
\end{figure}

\begin{figure}[hbtp]
\centering 
\includegraphics[width=4.5in]{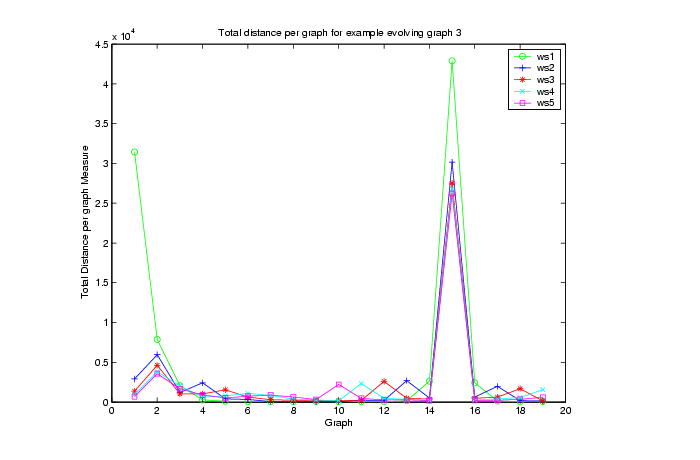}
\caption{The total distance per graph measure with different window sizes on evolving graph 3}
\label{fig:vs_2_ex3_g}
\end{figure}

\begin{figure}[hbtp]
\centering 
\includegraphics[width=4.5in]{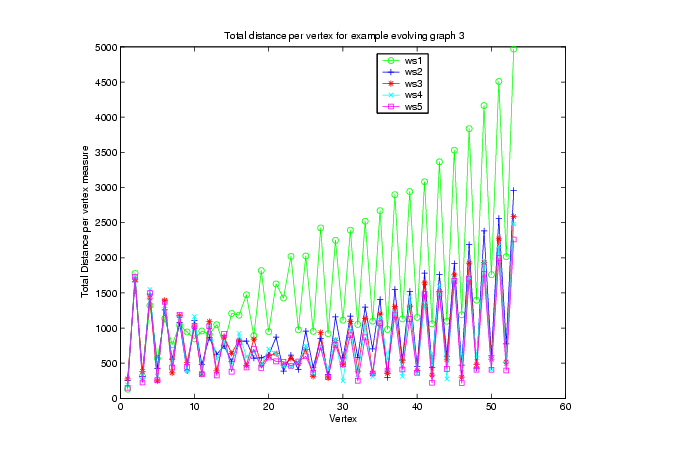}
\caption{The total distance per vertex measure with different window sizes on evolving graph 3}
\label{fig:vs_2_ex3_v}
\end{figure}

\begin{figure}[hbtp]
\centering 
\includegraphics[width=4.5in]{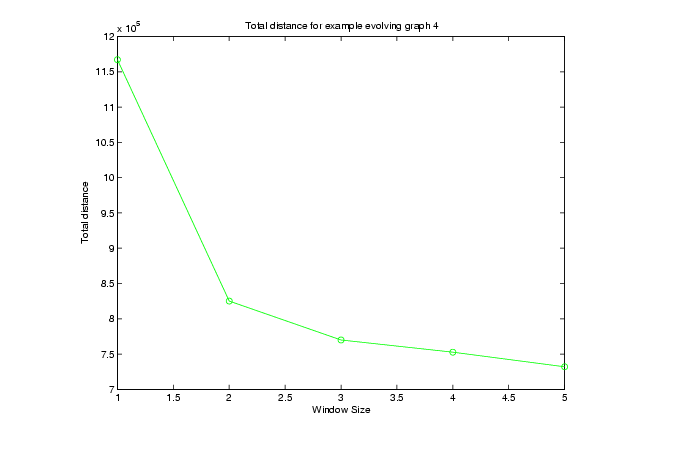}
\caption{The total distance measure with different window sizes on evolving graph 4}
\label{fig:vs_2_ex4_eg}
\end{figure}

\begin{figure}[hbtp]
\centering 
\includegraphics[width=4.5in]{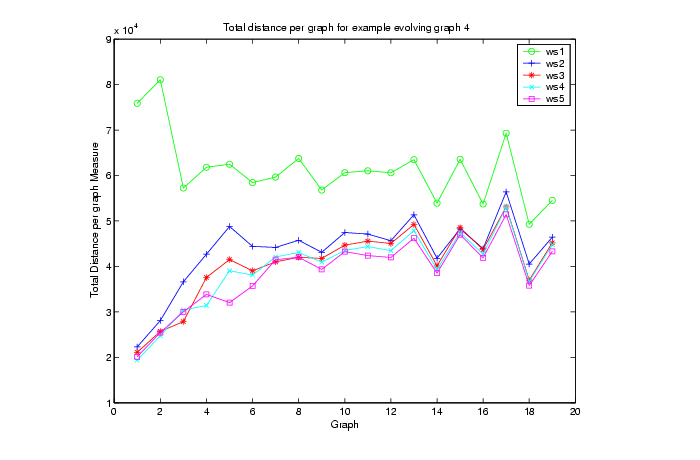}
\caption{The total distance per graph measure with different window sizes on evolving graph 4}
\label{fig:vs_2_ex4_g}
\end{figure}

\begin{figure}[hbtp]
\centering 
\includegraphics[width=4.5in]{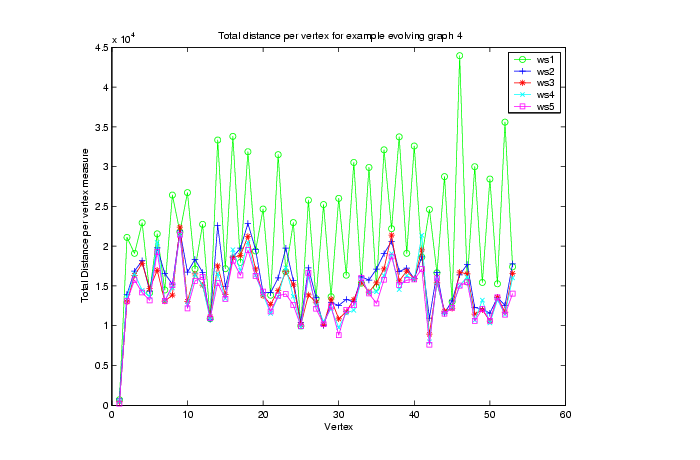}
\caption{The total distance per vertex measure with different window sizes on evolving graph 4}
\label{fig:vs_2_ex4_v}
\end{figure}

\begin{figure}[hbtp]
\centering 
\includegraphics[width=4.5in]{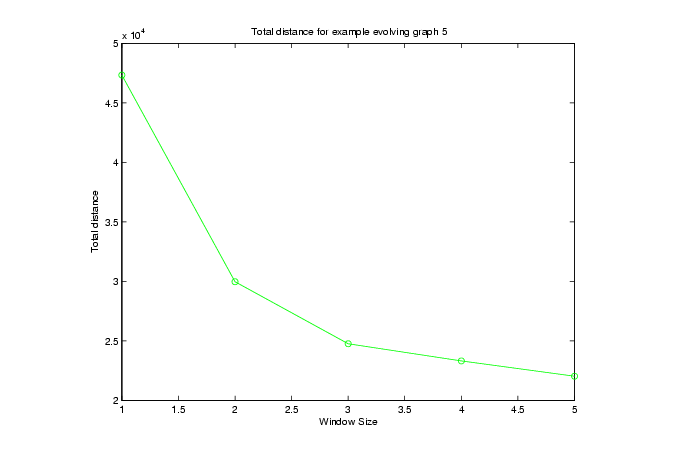}
\caption{The total distance measure with different window sizes on evolving graph 5}
\label{fig:vs_2_ex5_eg}
\end{figure}

\begin{figure}[hbtp]
\centering 
\includegraphics[width=4.5in]{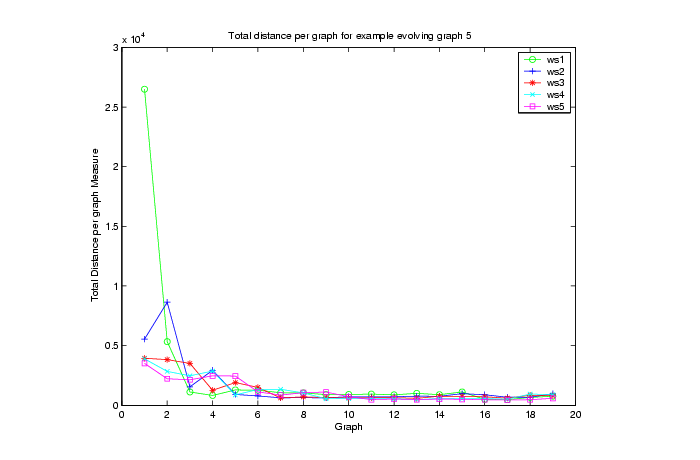}
\caption{The total distance per graph measure with different window sizes on evolving graph 5}
\label{fig:vs_2_ex5_g}
\end{figure}

\begin{figure}[hbtp]
\centering 
\includegraphics[width=4.5in]{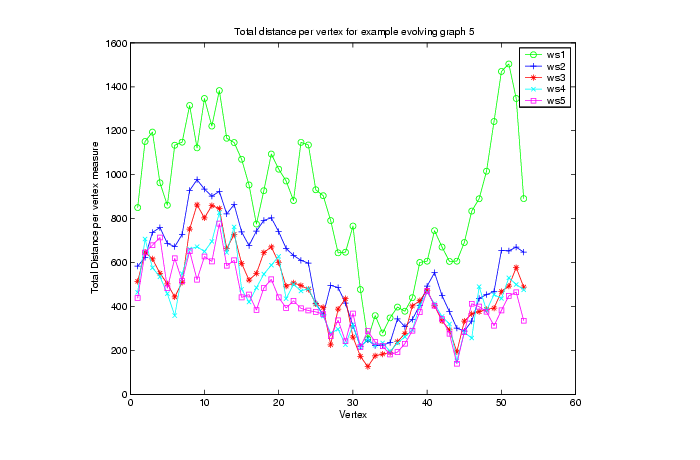}
\caption{The total distance per vertex measure with different window sizes on evolving graph 5}
\label{fig:vs_2_ex5_v}
\end{figure}

\begin{figure}[hbtp]
\centering 
\includegraphics[scale=0.1]{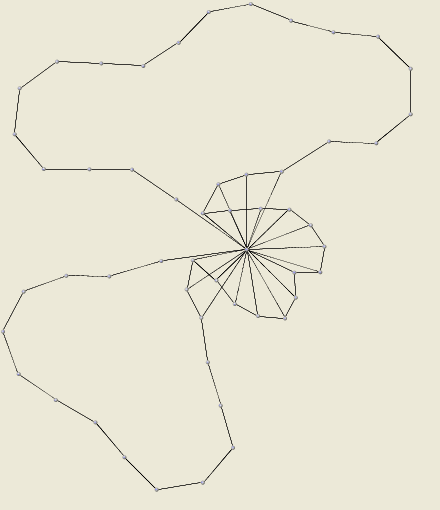}
\includegraphics[scale=0.1]{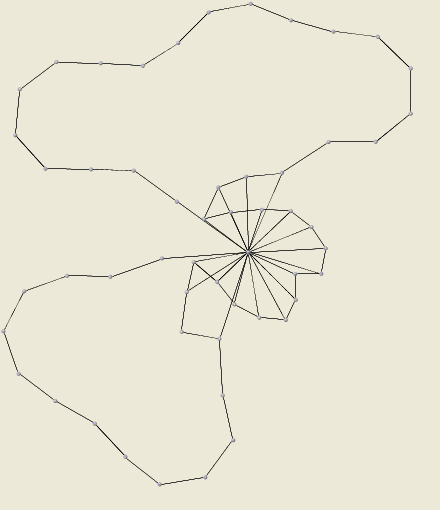}
\includegraphics[scale=0.1]{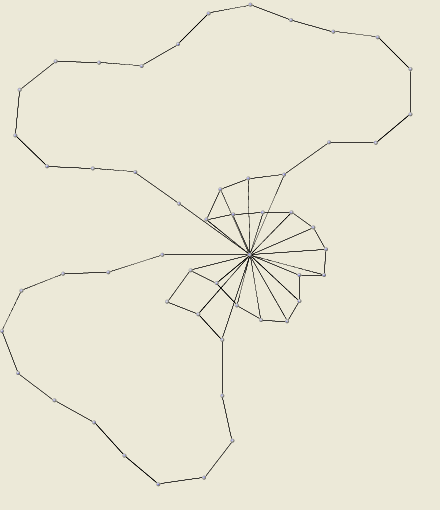}
\includegraphics[scale=0.1]{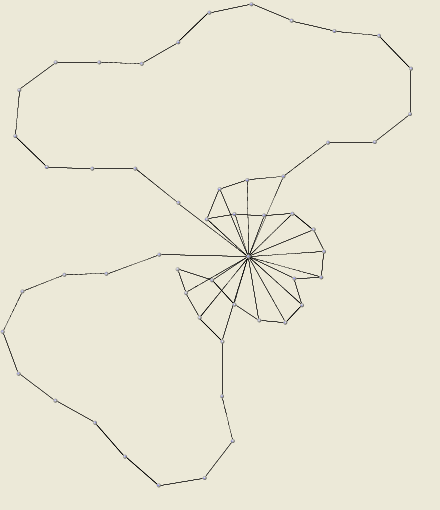}
\includegraphics[scale=0.1]{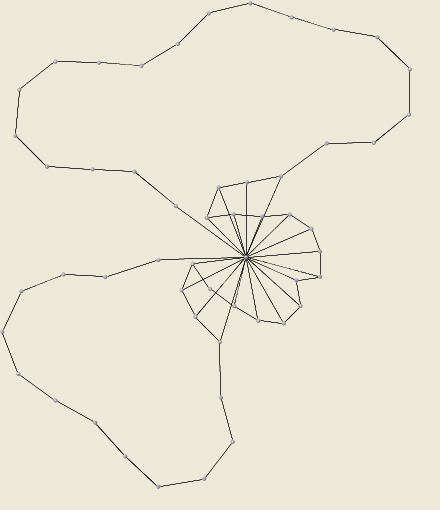}
\includegraphics[scale=0.1]{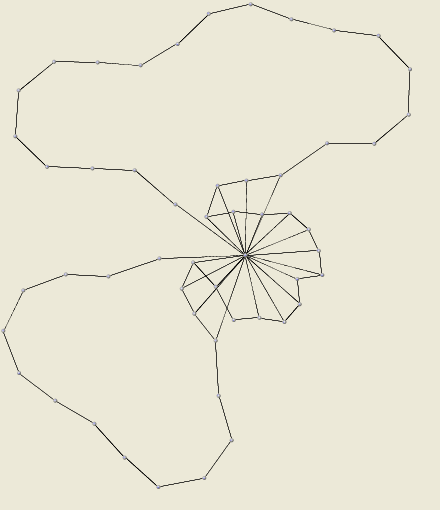}
\includegraphics[scale=0.1]{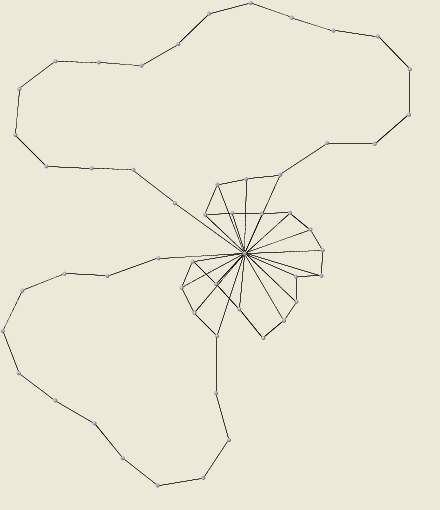}
\includegraphics[scale=0.1]{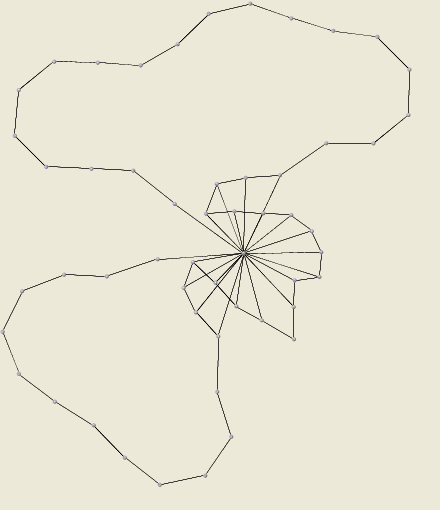}
\\
\includegraphics[scale=0.1]{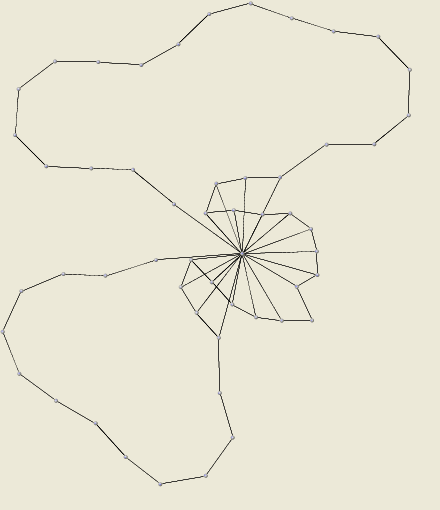}
\includegraphics[scale=0.1]{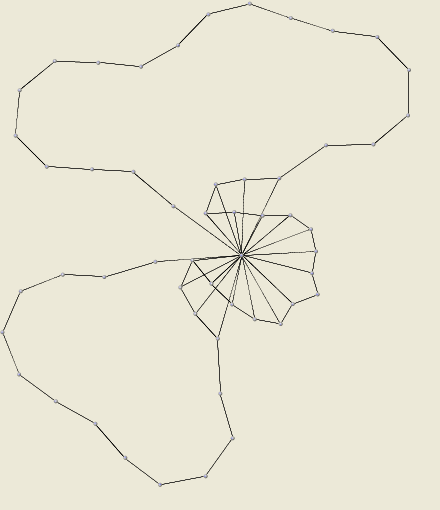}
\includegraphics[scale=0.1]{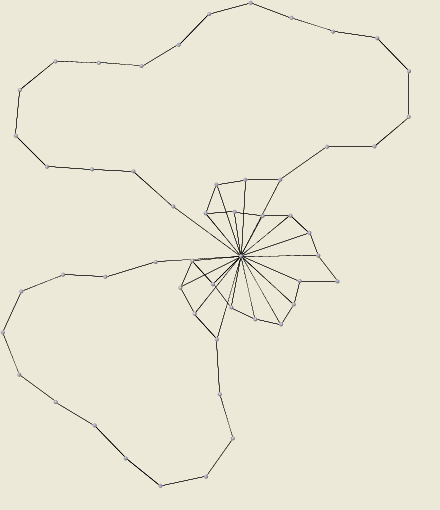}
\includegraphics[scale=0.1]{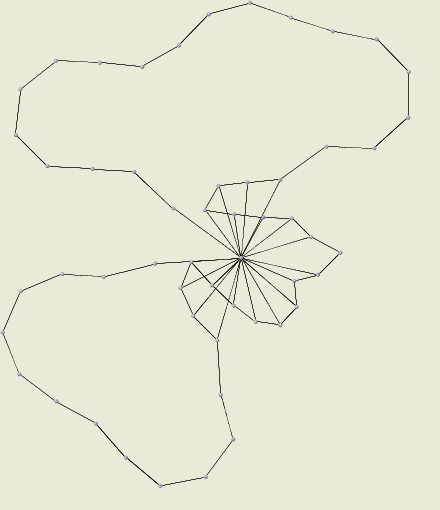}
\includegraphics[scale=0.1]{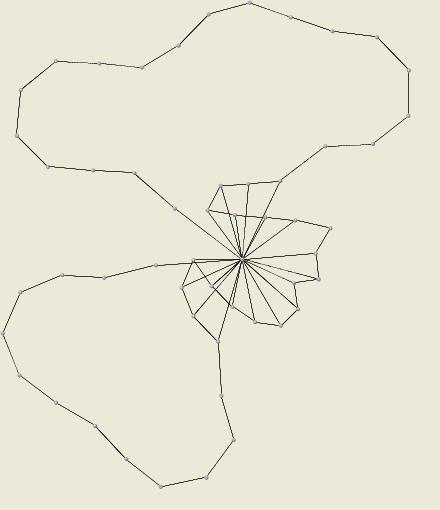}
\includegraphics[scale=0.1]{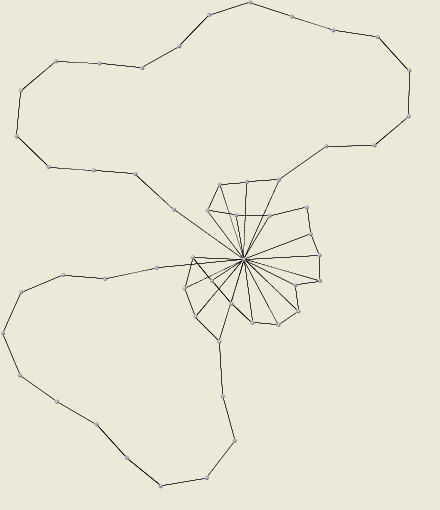}
\includegraphics[scale=0.1]{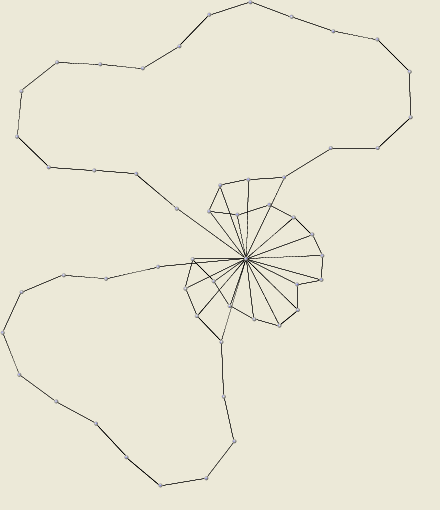}
\includegraphics[scale=0.1]{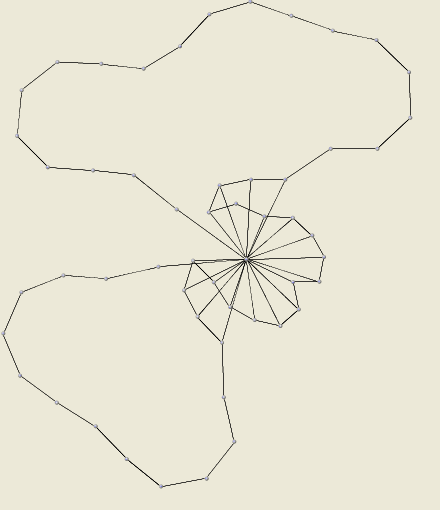}
\\
\includegraphics[scale=0.1]{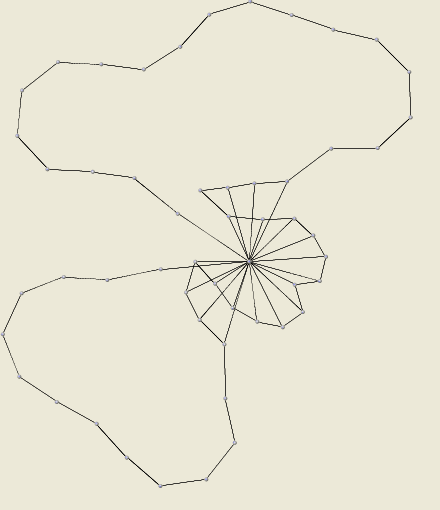}
\includegraphics[scale=0.1]{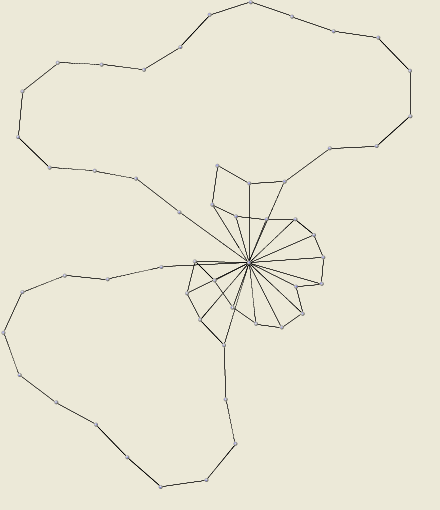}
\includegraphics[scale=0.1]{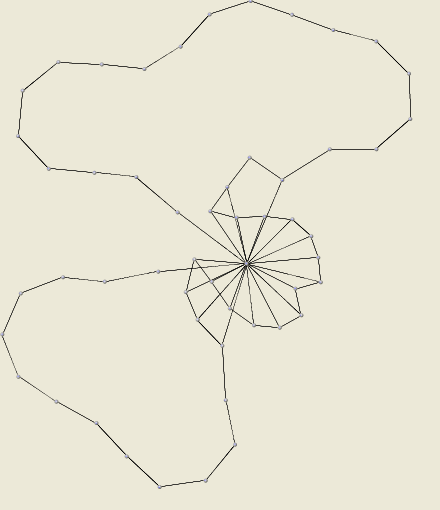}
\includegraphics[scale=0.1]{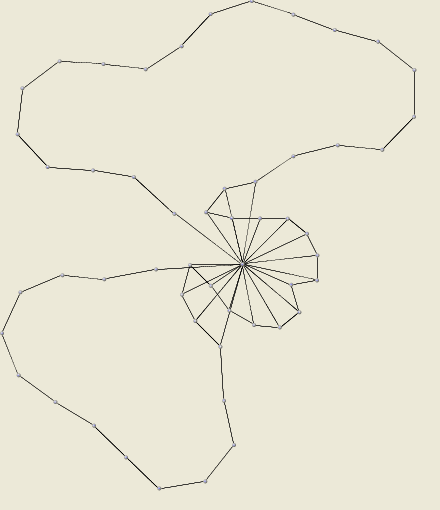}
\\
\caption{Example evolving graph 5 with vertex optimization layout}
\label{fig:app_eg_ex5_opt}
\end{figure}

{\bf Discussion}

Figures \ref{fig:vs_2_ex1_eg}, \ref{fig:vs_2_ex2_eg}, \ref{fig:vs_2_ex3_eg},
\ref{fig:vs_2_ex4_eg}, and \ref{fig:vs_2_ex5_eg} show the total distance measure with window size from one to five on evolving graphs 1 to 5, respectively.
Figures \ref{fig:vs_2_ex1_g}, \ref{fig:vs_2_ex2_g}, \ref{fig:vs_2_ex3_g},
\ref{fig:vs_2_ex4_g}, and \ref{fig:vs_2_ex5_g} show the total distance per graph measure with window size from one to five on evolving graphs 1 to 5, respectively.
Figures \ref{fig:vs_2_ex1_v}, \ref{fig:vs_2_ex2_v}, \ref{fig:vs_2_ex3_v},
\ref{fig:vs_2_ex4_v}, and \ref{fig:vs_2_ex5_v} show the total distance per vertex measure with window size from one to five on evolving graphs 1 to 5, respectively.

In Figure \ref{fig:vs_2_ex1_eg} the plot of total distance measure drops from
window size of 1 to 2 and is quite steady afterward. This suggests that window size of 1 is not big enough and after window size of 2, there is small gain on improvement.

Figure \ref{fig:vs_2_ex1_g} shows that most of the total distance per graph measure for window size of 1 comes from the difference between the first and second graph instance. All window sizes have high total distance per graph between the second and third graph instances. The figure shows that the first couple of graph instances are in unsteady state because the force directed algorithm moves vertices in the first graph instance further in the early round, but are unable to move them back since the first graph instance is constrained by other graph instances which are still not in their optimized position; contrary to the other graph instances which are constrained by some optimized graph instances, which, therefore, are likely to be in a more steady state.

Figure \ref{fig:vs_2_ex1_v} shows that total distance per vertex measure on evolving graph 1 with window size of 1 keeps increasing for later vertices where other window size have steady total distance measure for all vertices. This suggests that the vertices are moved away further at the end of a line graph. 

Figure \ref{fig:vs_2_ex2_eg} is a plot of total distance measure on evolving graph 2. It shows the drop as window size increases. However, the total distance measure is significantly higher than the plot of evolving graph 1.

In Figure \ref{fig:vs_2_ex2_g}, the plot of total distance measure per graph on evolving graph 2 with window size 1 shows higher total distance than the other window sizes for all graph instances. The plot also shows the same effect between graph instance 2 and 3 as in Figure \ref{fig:vs_2_ex1_g}.

In Figure \ref{fig:vs_2_ex2_v}, the plot of total distance measure per vertex on evolving graph 2 shows varying total distance between vertices. This shows that all the vertices are disturbed by the configuration of switching between two sets of graph instances.

In Figure \ref{fig:vs_2_ex3_eg}, the plot of total distance measure on evolving graph 3 shows the drop as window size increases. The total distance measure is higher than that of evolving graph 1.

In Figure \ref{fig:vs_2_ex3_g}, the plot of total distance measure per graph on evolving graph 3 shows the high total distance measure between the two graph instance with different configurations, which is expected. It is significant to point out that another high total distance occurs at another graph instance far from the highest total distance by the amount of window size. For example, with window size of 4, the highest total distance is at graph 15; graphs 11 and 19 also have high total distance. 

In Figure \ref{fig:vs_2_ex3_v}, the plot of total distance measure per vertex on evolving graph 3 has an interesting pattern. The odd vertices seem to have low total distance and the even vertices have high total distance. From the configuration, this is reasonable because the first set of graph instances connects the odd and even vertices together and the second set of graph instances connects the odd vertices with odd vertices and the even vertices with even vertices.

In Figure \ref{fig:vs_2_ex4_eg}, the plot of total distance measure on evolving graph 4 has the same shape as that of other evolving graphs. The total distance value, however, is significantly higher than the other evolving graphs. Since all graph instances have different configurations, this is expected.

In Figure \ref{fig:vs_2_ex4_g}, the plot of total distance measure per graph on evolving graph 4 shows high total distance for all graph instances for all window sizes. Window sizes 2 to 5 all have a comparatively similar plot. The configuration is set up so that the same vertex is farther for later graph instances. Therefore, the plot increases from graph instance 1 to 5. After graph instance 5, with the roundup, the vertices are no farther any more and total distance stops increasing.

In Figure \ref{fig:vs_2_ex4_v}, the plot of total distance measure per vertex on evolving graph 4 shows high total distance for window size of 1. Window sizes 2 to 5 have comparatively similar plots.

In Figure \ref{fig:vs_2_ex5_eg}, the plot of total distance measure on evolving graph 5 drops as the window size increases.

In Figure \ref{fig:vs_2_ex5_g}, the plot of total distance measure per graph on evolving graph 5 shows unsteady state in the first couple of graph instances and steady state afterward for all window sizes. The difference of one edge seems to have little impact on the overall placement.

In Figure \ref{fig:vs_2_ex5_v}, the plot of total distance measure per vertex on evolving graph 5 shows varying total distance for each vertex. Vertex 40 with the most changes has high total distance but does not have the highest total distance. Vertices 1 to 20 which connect to vertex 40 all have high total distances. The plot around vertices 1 to 20 looks like a bell curve and the plot around vertex 40 looks like a spike. This suggests that vertices around the high total distance vertices are pulled. This is reasonable because the graph configuration is circular. The graph layout using vertex optimization algorithm for evolving graph 5 is in Figure \ref{fig:app_eg_ex5_opt}.

Overall, the results show that for simple evolving graphs such as example 1 and example 2, where the best result is to retain the position of all vertices, the vertex optimization algorithm might not obtain the best result. However, considering other factors such as the optimum position of vertices relative to other vertices in the same graph instances that the force directed algorithm provides, the vertex optimization algorithm can be used for automatic evolving graph layout.

Different window sizes yield different results with high value of window size generally resulting in better total distance measure. Window size 1 gives poor results compared to others and window size 2 is adequate for very simple evolving graphs.

Evolving graph 1 shows that error in the first few graph instances contributes to the total distance measure. Evolving graph 2 shows that the difference in graph configuration contributes to the total distance measure. Evolving graphs 3, 4 and 5 show that the total distance per graph measure can tell roughly where the high change in graph configuration occurs. Total distance per vertex measure can explain the difference in graph configuration.

\section*{Appendix B}

\subsection*{EGML DTD Specifications}

\begin{verbatim}
<!-- DTD for EGML 1.0 -->
<!-- Authors: Anurat Chapanond -->
<!-- Computer Science Department -->
<!-- Rensselaer Polytechnic Institute -->
<!-- egml.dtd,v 1.0 03/16/2007 -->

<!-- Boolean type -->
<!ENTITY % boolean.type "(0|1)" >

<!-- Positive number type -->
<!ENTITY % number.type "NMTOKEN">

<!-- ID type -->
<!ENTITY % id.type "NMTOKEN">

<!-- String type -->
<!ENTITY % string.type "CDATA">


<!-- Standard XML Namespace attribute -->
<!ENTITY % nds 'xmlns'>

<!-- URI type -->
<!ENTITY % uri.type "%string.type;">

<!-- Anchor type -->
<!ENTITY % anchor.type "(c|n|ne|e|se|s|sw|w|nw)">

<!-- Type of Graphics (GML types) -->
<!ENTITY % type-graphics-gml.type "arc|bitmap|image|line|
				oval|polygon|rectangle|text">

<!-- Type of Graphics (New types) -->
<!ENTITY % type-graphics-app.type "box|circle|ver_ellipsis|
				hor_ellipsis|rhombus|triangle|pentagon|
				hexagon|octagon">

<!-- Line types -->
<!-- Arrow type -->
<!ENTITY % arrow.type "(none | first | last | both)">
<!-- Capstyle type -->
<!ENTITY % capstyle.type "(butt | projecting | round)">
<!-- Joinstyle type -->
<!ENTITY % joinstyle.type "(bevel | miter | round)">
<!-- Arc style  type -->
<!ENTITY % arcstyle.type "(pieslice | chord | arc)">

<!-- Text types -->
<!-- Text justification type -->
<!ENTITY % justify.type "(left | right | center)">
<!-- Font type -->
<!ENTITY % font.type "%string.type;">
<!-- Color type -->
<!ENTITY % color.type "%string.type;">

<!-- Angle type -->
<!ENTITY % angle.type "%string.type;">

<!-- Object type -->
<!ENTITY % object.type "(list | string | real | integer)">

<!-- Global Attributes -->
<!ENTITY % global-atts 
		"id %number.type; #IMPLIED
		name %string.type; #IMPLIED
		label %string.type; #IMPLIED
		labelanchor %string.type; #IMPLIED">

<!-- Standard XML Attributes -->
<!ENTITY % xml-atts "%nds; %uri.type; #FIXED 
					'http://www.cs.rpi.edu/XGMML' 
                    xml:lang NMTOKEN #IMPLIED
                    xml:space (default | preserve) #IMPLIED">

<!-- Standard XLink Attributes -->
<!ENTITY % xlink-atts 
              "xmlns:xlink CDATA #FIXED 'http://www.w3.org/1999/xlink'
               xlink:type (simple) #FIXED 'simple' 
               xlink:role CDATA #IMPLIED
               xlink:title CDATA #IMPLIED
               xlink:show (new|embed|replace) #FIXED 'replace'
               xlink:actuate (onLoad|onRequest) #FIXED 'onRequest'
               xlink:href CDATA #IMPLIED">

<!-- Safe Graph Attributes -->

<!ENTITY % graph-atts-safe "directed %boolean.type; '0' ">		
<!-- Unsafe Graph Attributes (GML) -->
<!ENTITY % graph-atts-gml-unsafe "Vendor %string.type;  #IMPLIED">

<!-- Unsafe Graph Attributes (new attributes) -->
<!ENTITY % graph-atts-app-unsafe "Scale %number.type; #IMPLIED
                                  Rootnode %number.type; #IMPLIED
                                  Layout %string.type; #IMPLIED
                                  Graphic %boolean.type; #IMPLIED">

<!-- Graph Element -->
<!ELEMENT graph (att*,(node | edge)*)>
<!-- Graph Attributes -->
<!ATTLIST graph
	%global-atts;
	%xml-atts;
	%xlink-atts;
	%graph-atts-safe;
	%graph-atts-gml-unsafe;
	%graph-atts-app-unsafe;>

<!-- Safe Node Attributes (GML) -->
<!ENTITY % node-atts-gml-safe "edgeanchor %string.type; #IMPLIED">
<!-- Safe Node Attributes (new attributes) -->
<!ENTITY % node-atts-app-safe "weight %string.type;  #IMPLIED">

<!-- Node Element -->
<!ELEMENT node (graphics?,att*)>
<!-- Node Attributes -->
<!ATTLIST node
	%global-atts;
	%xlink-atts;
	%node-atts-gml-safe;
	%node-atts-app-safe;>

<!-- Safe Edge Attributes (GML) -->
<!ENTITY % edge-atts-gml-safe "source %number.type; #REQUIRED
                               target %number.type; #REQUIRED">
<!-- Safe Edge Attributes (new attributes) -->
<!ENTITY % edge-atts-app-safe "weight %string.type; #IMPLIED">

<!-- Edge Element -->
<!ELEMENT edge (graphics?,att*)>
<!-- Edge Attributes -->
<!ATTLIST edge
	%global-atts;
	%xlink-atts;
	%edge-atts-gml-safe;
	%edge-atts-app-safe;>

<!-- Graphics Type  -->
<!ENTITY % graphics-type-att "type (%type-graphics-gml.type;|
				%type-graphics-app.type;) #IMPLIED">

<!-- Point Attributes (x,y,z)  -->
<!ENTITY % point-atts "x %number.type; #IMPLIED
				y %number.type; #IMPLIED
				z %number.type; #IMPLIED">

<!-- Dimension Attributes (width,height,depth)  -->
<!ENTITY % dimension-atts "w %number.type; #IMPLIED
				h %number.type; #IMPLIED
				d %number.type; #IMPLIED">

<!-- External Attributes (Image and Bitmap)  -->
<!ENTITY % external-atts "image %uri.type; #IMPLIED
				bitmap %uri.type; #IMPLIED">

<!-- Line Attributes -->
<!ENTITY % line-atts "width %number.type; #IMPLIED
				arrow %arrow.type; #IMPLIED
				capstyle %capstyle.type; #IMPLIED
				joinstyle %joinstyle.type; #IMPLIED
				smooth %boolean.type; #IMPLIED
				splinesteps %number.type; #IMPLIED">

<!-- Text Attributes -->
<!ENTITY % text-atts "justify %justify.type; #IMPLIED
				font  %font.type; #IMPLIED">

<!-- Bitmap Attributes -->

<!ENTITY % bitmap-atts "background %color.type; #IMPLIED
				foreground %color.type; #IMPLIED">

<!-- Arc Attributes -->
<!ENTITY % arc-atts "extent %angle.type; #IMPLIED
				start %angle.type; #IMPLIED
				style %arcstyle.type;  #IMPLIED">

<!-- Graphical Object Attributes -->
<!ENTITY % object-atts "stipple %string.type; #IMPLIED
				visible %boolean.type; #IMPLIED
				fill %color.type; #IMPLIED
				outline %color.type; #IMPLIED
				anchor %anchor.type; #IMPLIED">

<!-- Graphics Element -->
<!ELEMENT graphics ((Line? | center?),att*)>
<!-- Graphics Attributes -->
<!ATTLIST graphics
          %graphics-type-att;
	  %point-atts;
	  %dimension-atts;
	  %external-atts;
	  %line-atts;
	  %text-atts;
	  %bitmap-atts;
	  %arc-atts;
	  %object-atts;>

<!-- Center Point Element -->
<!ELEMENT center EMPTY>
<!ATTLIST center 
	  %point-atts;>

<!-- Line Element -->
<!ELEMENT Line (point,point+)>

<!-- Point Element -->
<!ELEMENT point EMPTY>
<!ATTLIST point 
	  %point-atts;>

<!-- Value Attribute -->
<!ENTITY % attribute-value "value %string.type; #IMPLIED">
<!-- Type Attribute -->
<!ENTITY % attribute-type "type %object.type;  #IMPLIED">

<!-- Att Element -->
<!ELEMENT att (#PCDATA | att | graph)*>
<!-- Att Attributes -->
<!ATTLIST att
          %global-atts;
	  %attribute-value;
	  %attribute-type;>

<!-- Algorithm Attributes -->
<!ENTITY % algo-atts "algorithm-name %string.type; #IMPLIED">

<!-- Counting the number Attributes -->
<!ENTITY % count-atts "no	%number.type; #IMPLIED">

<!ELEMENT	evolving-graph (graph-instance*,prediction?)>
<!ATTLIST	evolving-graph	
	%global-atts;
	%count-atts;
>

<!ELEMENT	graph-instance	(graph,timestamp,metric*,cluster?,rank?)>
<!ATTLIST	graph-instance	%global-atts;>

<!ELEMENT	timestamp	(#PCDATA)>

<!ELEMENT	metric	(value|value-list)>
<!ATTLIST	metric
	%global-atts;
	%algo-atts;
>

<!ELEMENT	value	(#PCDATA)>
<!ELEMENT	value-list	(value*)>
<!ATTLIST	value-list
	%global-atts;
	%count-atts;
>

<!ELEMENT	cluster	(node-set*)>
<!ATTLIST	cluster
	%global-atts;
	%algo-atts;
	%count-atts;
>

<!ELEMENT	node-set	(node*)>
<!ATTLIST	node-set
	%global-atts;
	%count-atts;
>


<!ELEMENT	rank	(rank-element*)>
<!ATTLIST	rank
	%global-atts;
	%algo-atts;
>

<!ELEMENT	rank-element	((node|edge),position,value)>
<!ELEMENT	position	(#PCDATA)>

<!ELEMENT	prediction	(timestamp,rank-element*)>
<!ATTLIST	prediction
	%global-atts;
	%algo-atts;
>
\end{verbatim}

\subsection*{EGML Schema Specifications}

\begin{verbatim}
<?xml version='1.0' encoding='UTF-8'?>

<!-- XML schema for EGML 1.0 -->
<!-- Authors: Anurat Chapanond -->
<!-- Computer Science Department -->
<!-- Rensselaer Polytechnic Institute -->
<!-- egml.xsd,v 1.0 03/20/2007 -->

<xsd:schema xmlns:xsd="http://www.w3.org/2001/XMLSchema" 
	targetNamespace="http://www.cs.rpi.edu/XGMML" 
	xmlns="http://www.cs.rpi.edu/XGMML" 
	xmlns:xml="http://www.w3.org/XML/1998/namespace"
	xmlns:xlink="http://www.w3.org/1999/xlink"
	elementFormDefault="qualified"
	attributeFormDefault="unqualified"
	version="egml 1.0">

<!-- Include XGMML schema -->
<xsd:include schemaLocation="http://www.cs.rpi.edu/~puninj/XGMML/xgmml.xsd" />

<!-- Algorithm Name Attributes -->
<xsd:attributeGroup name="algo-atts" >
	<xsd:attribute name="algorithm-name" type="string.type" />
</xsd:attributeGroup>

<!-- Counting the number Attributes -->
<xsd:attributeGroup name="count-atts" >
	<xsd:attribute name="no" type="number.type" />
</xsd:attributeGroup>

<!-- Evolving graph element -->
<xsd:element name="evolving-graph" type="evolving-graph.type" />
<xsd:complexType name="evolving-graph.type" >
	<xsd:sequence>
		<xsd:element ref="graph-instance" minOccurs="0" maxOccurs="unbounded" />
		<xsd:element ref="prediction" minOccurs="0" maxOccurs="1" />
	</xsd:sequence>
	<xsd:attributeGroup ref="global-atts" />
	<xsd:attributeGroup ref="count-atts" />
</xsd:complexType>

<!-- Graph instance element -->
<xsd:element name="graph-instance" type="graph-instance.type" />
<xsd:complexType name="graph-instance.type" >
	<xsd:sequence>
		<xsd:element ref="graph" minOccurs="1" maxOccurs="1" />
		<xsd:element ref="timestamp" minOccurs="1" maxOccurs="1" />
		<xsd:element ref="metric" minOccurs="0" maxOccurs="unbounded" />
		<xsd:element ref="cluster" minOccurs="0" maxOccurs="1" />
		<xsd:element ref="rank" minOccurs="0" maxOccurs="1" />
	</xsd:sequence>
	<xsd:attributeGroup ref="global-atts" />
</xsd:complexType>

<!-- Timestamp element -->
<xsd:element name="timestamp" type="number.type" />

<!-- Metric element -->
<xsd:element name="metric" type="metric.type" />
<xsd:complexType name="metric.type" >
	<xsd:choice>
		<xsd:element ref="value" minOccurs="0" maxOccurs="1" />
		<xsd:element ref="value-list" minOccurs="0" maxOccurs="1" />
	</xsd:choice>
	<xsd:attributeGroup ref="global-atts" />
	<xsd:attributeGroup ref="algo-atts" />
</xsd:complexType>

<!-- Value element -->
<xsd:element name="value" type="xsd:double" />

<!-- Value-list element -->
<xsd:element name="value-list" type="value-list.type" />
<xsd:complexType name="value-list.type" >
	<xsd:sequence>
		<xsd:element ref="value" minOccurs="0" maxOccurs="unbounded" />
	</xsd:sequence>
	<xsd:attributeGroup ref="global-atts" />
	<xsd:attributeGroup ref="count-atts" />
</xsd:complexType>

<!-- Cluster element -->
<xsd:element name="cluster" type="cluster.type" />
<xsd:complexType name="cluster.type" >
	<xsd:sequence>
		<xsd:element ref="node-set" minOccurs="0" maxOccurs="unbounded" />
	</xsd:sequence>
	<xsd:attributeGroup ref="global-atts" />
	<xsd:attributeGroup ref="algo-atts" />
	<xsd:attributeGroup ref="count-atts" />
</xsd:complexType>

<!-- Node-set element -->
<xsd:element name="node-set" type="node-set.type" />
<xsd:complexType name="node-set.type" >
	<xsd:sequence>
		<xsd:element ref="node" minOccurs="0" maxOccurs="unbounded" />
	</xsd:sequence>
	<xsd:attributeGroup ref="global-atts" />
	<xsd:attributeGroup ref="count-atts" />
</xsd:complexType>

<!-- Rank element -->
<xsd:element name="rank" type="rank.type" />
<xsd:complexType name="rank.type" >
	<xsd:sequence>
		<xsd:element ref="rank-element" minOccurs="0" maxOccurs="unbounded" />
	</xsd:sequence>
	<xsd:attributeGroup ref="global-atts" />
	<xsd:attributeGroup ref="algo-atts" />
</xsd:complexType>

<!-- Rank-element element -->
<xsd:element name="rank-element" type="rank-element.type" />
<xsd:complexType name="rank-element.type" >
	<xsd:sequence>
		<xsd:choice>
			<xsd:element ref="node" minOccurs="1" maxOccurs="1" />
			<xsd:element ref="edge" minOccurs="1" maxOccurs="1" />
		</xsd:choice>
		<xsd:element ref="position" minOccurs="1" maxOccurs="1" />
		<xsd:element ref="value" minOccurs="1" maxOccurs="1" />
	</xsd:sequence>
</xsd:complexType>

<!-- Position element -->
<xsd:element name="position" type="number.type" />

<!-- Prediction element -->
<xsd:element name="prediction" type="prediction.type" />
<xsd:complexType name="prediction.type" >
	<xsd:sequence>
		<xsd:element ref="timestamp" minOccurs="1" maxOccurs="1" />
		<xsd:element ref="rank-element" minOccurs="0" maxOccurs="unbounded" />
	</xsd:sequence>
	<xsd:attributeGroup ref="global-atts" />
	<xsd:attributeGroup ref="algo-atts" />
</xsd:complexType>

</xsd:schema>
\end{verbatim}

\bibliographystyle{acm}
\bibliography{../thesis/anuratBib}

\begin{thebibliography}{10}

\bibitem{Batagelj}
{\sc Batagelj, V., and Mrvar, A.}
\newblock Pajek – program for large network analysis.
\newblock {\em Connections 21(2)\/}, 47--57.

\bibitem{Bearman}
{\sc Bearman, P.~S., and Everett, K.}
\newblock The structure of social protest.
\newblock {\em Social Networks 15\/} (1993), 171--200.

\bibitem{Bender-deMoll}
{\sc Bender-deMoll, S., and McFarland, D.~A.}
\newblock Sonia(social network image animator).
\newblock {\em See http://sonia.stanford.edu\/} (2002).

\bibitem{Bradley}
{\sc Bradley, G.~H.}
\newblock{Network and Graph Markup Language File Formats}
\newblock {\em Nineth INFORMS Computing Society Conference Annapolis, Maryland, 2005}

\bibitem{Cole}
{\sc Cole, G., Bulashova, N., and Yurcik, W.}
\newblock Geographical netflows visualization for network situational
  awareness: Naukanet administrative data analysis system (nadas).
\newblock {\em 12th International Conference on Telecommunication Systems -
  Modeling and Analysis (ICTSM)\/} (2004).

\bibitem{Eurovision}
{\sc Eurovision}.
\newblock Eurovision song contest.
\newblock
{  http://web.archive.org/web/20060525094524/

http://www.eurovision.tv/english/2513.htm.
}
\bibitem{Fenn}
{\sc Fenn, D., Suleman, O., Efstathiou, J., and Johnson, N.~F.}
\newblock How does europe make its mind up? connections cliques and
  compatibility between countries in the eurovision song contest.

\bibitem{Fruchterman}
{\sc Fruchterman, T. M.~J., and Reingold, E.~M.}
\newblock Graph drawing by force-directed placement.
\newblock Tech. Rep. UIUCDCS-R-90-1609, Department of Computer Science,
  University of Illinois, Urbana-Champagne, Ill, June 1990.

\bibitem{Rollcall}
{\sc LibraryOfCongress}.
\newblock Roll call votes for us house of representatives and senate.
\newblock http://thomas.loc.gov/home/rollcallvotes.html.

\bibitem{Moody}
{\sc Moody, J., McFarland, D., and Bender-deMoll, S.}
\newblock Dynamic network visualization.
\newblock {\em American Journal of Sociology 110\/} (2005), 1206--1241.

\bibitem{Preston}
{\sc Preston, N., and Krishnamoorthy, M.~S.}
\newblock Graphdraw- a graph drawing system to study social networks.
\newblock GraphDraw Java based system, 2004.

\bibitem{PuninPhD}
{\sc Punin, J.~R.}
\newblock WWWPAL suite for analysis and organization of web sites. 
\newblock Doctoral Dissertation, RPI, Troy, NY, August 2003.

\bibitem{Punin3}
{\sc Punin, J.~R.}
\newblock XGMML (extensible graph markup and modeling language).
\newblock http://www.cs.rpi.edu/~puninj/XGMML/.

\bibitem{Punin}
{\sc Punin, J.~R., Krishnamoorthy, M.~S., and Zaki, M.}
\newblock Web usage mining – languages and algorithms.
\newblock {\em WebKDD Workshop, ACM SIGKDD\/} (2001), 88--112.

\bibitem{Tollis}
{\sc Tollis, I.~G., Battista, G., Eades, P., and Tamassia, R.}
\newblock {\em Graph Drawing: Algorithms for the Visualization of Graphs},
  1st~ed.
\newblock Prentice Hall, 1998.

\bibitem{Watts}
{\sc Watts, D.~J., Dodds, P.~S., and Newman, M. E.~J.}
\newblock Identity and search in social networks.
\newblock {\em http://arxiv.org/abs/cond-mat/0205383\/} (May 2006).

\bibitem{Wegman}
{\sc Wegman, E., and Marchette, D.}
\newblock On some techniques for streaming data: A case study of internet
  packet headers.
\newblock {\em JCGS 12(4)\/}, 893--914.

\bibitem{XGMML}
\newblock{ \em http://en.wikipedia.org/wiki/XGMML\/} 
\end{thebibliography}

\end{document}